\documentclass[a4paper]{article}
\usepackage{geometry}
\usepackage[T1]{fontenc}
\usepackage[utf8]{inputenc}
\usepackage[english]{babel}
\usepackage{mathtools}
\usepackage{amsmath}
\usepackage{amssymb}
\usepackage{xcolor}
\usepackage{graphicx}
\usepackage{nicefrac}
\usepackage[backend=bibtex,sorting=none]{biblatex}
\usepackage{color, colortbl}
\usepackage{algorithm,algorithmic}
\usepackage{hyperref} %
\usepackage{caption}
\usepackage{placeins}
\usepackage{tikz}
\usepackage{subcaption} 
\usepackage{pgfplots}
\usepackage{datetime}
\usepackage{tikz}

\newcommand{\imagi}{\mathrm{i}} %
\newcommand{\euler}{\mathrm{e}} %
\newcommand{\dxdv}{\,\mathrm{d}x \mathrm{d}v } %
\newcommand{\dx}{\,\mathrm{d}x }

\floatname{algorithm}{Algorithm}

\newcommand{\kvec}[1]{\ensuremath{\mathrm{vec}(#1)}}

\author{Jakob Ameres}

\author{
\vspace{.5em}\\
\large{Jakob Ameres$^{1,2}$} \\
\vspace{.5em}\\
$^1$\normalsize{Technische Universit\"at M\"unchen, Zentrum Mathematik}\\
\normalsize{Boltzmannstra\ss{}e 3, 85748 Garching, Deutschland}%
\vspace{.5em}\\
$^2$\normalsize{Max-Planck-Institut f\"ur Plasmaphysik}\\
\normalsize{Boltzmannstra\ss{}e 2, 85748 Garching, Deutschland}%
\vspace{.5em}\\
}

\title{Quasi Monte Carlo inverse transform sampling for phase space conserving Lagrangian particle methods and Eulerian-Lagrangian coupling}
\begin{document}
\maketitle 

\begin{abstract}
This article presents a novel and practically useful link between geometric integration, low-discrepancy sampling and code coupling for Lagrangian and Eulerian Vlasov-Poisson solvers.\\
Low-discrepancy sequences, also called quasi-random sequences (Quasi Monte Carlo), provide convergence rates close to $\mathcal{O}( N^{-1})$ which are far superior to (pseudo) random numbers (Monte Carlo) settling in at only $\mathcal{O}(N^{-0.5})$.
Lagrangian particle methods such as PIC rely on Monte Carlo integration. The particle distributions are nonlinearly perturbed by the forward flow following the characteristics.
Hence it remains the question of whether particle methods can benefit from such quasi-random-sequences. Any nonlinear measure-preserving map does not affect the low-discrepancy of a QMC sequence such that the order of convergence remains.
This article shows that the forward flow of phase space-conserving geometric particle methods induces naturally such a measure-preserving map underlying their importance in a new framework. In this context the Hardy Krause Variation is observed to increase in the Vlasov-Poisson system for the first time. with the linear phase. Also the star discrepancy is presented for an entire PIC simulation.\\
On the other hand, Eulerian and Lagrangian solvers have different strengths and weaknesses, such that we present a novel way of transiting from a spectral discretization of the Vlasov--Poisson system to a PIC simulation. This is achieved by higher dimensional inverse transform sampling (Rosenblatt-Mück transform). In this way Markov Chain Monte Carlo techniques are circumvented which allows the use of pseudo and quasi-random numbers. In the latter case better convergence rates can be observed both in the linear and nonlinear phase. \\ 
\end{abstract}


\newpage

\tableofcontents

\section{Introduction}
For the solution of kinetic models such as the two dimensional Vlasov--Poisson system~\eqref{vlasov1}-\eqref{poisson1}, Eulerian solvers and Lagrangian particle methods dominate the landscape.
\begin{align}
 \partial_t f(x,v,t) + v \cdot \partial_x f(x,v,t)  + \frac{q}{m} E(x,t) \cdot \partial_v f(x,v,t)=0 \label{vlasov1}\\
 E(x,t)=- \partial_x \Phi(x,t) \\
  \partial_{xx} \Phi(x,t) = q\Big(\int_{-\infty}^{\infty} f(x,v,t)  \mathrm{d}v -1 \Big) \label{poisson1}
\end{align}
Eulerian solvers represent the phase space density $f$ on a fixed grid. Here, we consider the simplest Eulerian solver, which is based on a Fourier spectral discretization 
of the entire phase space~\cite{joyce1971numerical,izrar1989integration} and a Hamiltonian splitting in time~\cite{watanabe2005vlasov}. Using the same time discretization~\cite{forest1989fourth} 
but discretizing the phase space by Monte Carlo samples one can obtain the standard geometric PIC (particle in cell)\cite{LEWIS1970136,evstatiev2013variational,kraus2013variational} method.
Unfortunately, as a Monte Carlo method PIC suffers from noise~\cite{ameres2018} with the slow $N^{-\frac{1}{2}}$ convergence, which can be improved to 
$N^{-(1-\epsilon) }$ by the use of Quasi Monte Carlo numbers~\cite{hickernell2005}. Still, the noise is especially a problem for small amplitudes in the initial phase of a simulation~\cite{ameres2018}.
Initially, the spectral solver appears to have no issues but after some time into the simulation, it suffers from the recurrence phenomenon \cite{einkemmer2014strategy} or 
filamentations~\cite{klimas1994splitting}, which can be mitigated by anti-aliasing and Fourier filtering techniques. 
Representing the density grid can be expensive in high dimensions~\cite{kormann2016sparse} and is also wasteful if large portions of phase space are practically empty. Here PIC performs better
as the markers can be placed with respect to the actual density which is known as importance sampling~\cite{caflisch98}.\\
\section{Particle in Cell}
Let us recall the fundamentals of the Particle in Cell (PIC) method~\cite{kraus2016gempic,birdsall2004plasma,decyk2011description,esirkepov2001exact}. Equation \eqref{vlasov1} describes a conservation law, which PIC solves by the methods of characteristics. 
\subsection{Method of characteristics}
The characteristics $\left( \mathbf{V}(t), \mathbf{X}(t) \right)$ are defined as a curve in space time along which the value of the density $f$ stays constant:
\begin{multline}
\label{densityf:char}
 \frac{\mathrm{d}}{\mathrm{d}t} f\left(\mathbf{X}(t),\mathbf{V}(t),t\right) = 
 \frac{\mathrm{d} \mathbf{X}(t) }{\mathrm{d}t}  \partial_x f \left(\mathbf{X}(t),\mathbf{V}(t),t \right)\\
 +\frac{\mathrm{d} \mathbf{V}(t) }{\mathrm{d}t}  \partial_v f \left(\mathbf{X}(t),\mathbf{V}(t),t \right)
 +\partial_t f \left(\mathbf{X}(t),\mathbf{V}(t),t \right) =0.
\end{multline}
Inserting $\partial_t f$ from~\eqref{vlasov1} into~\eqref{densityf:char} yields the equations
of motions for the characteristics of eqn.~\eqref{vlasov1}, which read
\begin{equation}
\label{vlasov1:char}
 \frac{\mathrm{d}}{\mathrm{d}t} \mathbf{V}(t) = - \frac{q}{m} E(t,\mathbf{X}(t)) 
 \text{ and } \frac{\mathrm{d}}{\mathrm{d}t} \mathbf{X}(t) = \mathbf{V}(t).
\end{equation}
Then $f$ as solution of eqn.~\eqref{vlasov1}
is constant along the characteristics~\eqref{vlasov1:char}, which means for given initial position in phase space
$(\mathbf{X}_0,\mathbf{V}_0)$ we have
\begin{equation}
\label{vlasov:characteristics}
 f(X(t=0),V(t=0),t=0) = f(\mathbf{X}(t),\mathbf{V}(t),t) \quad \forall t \geq 0.
\end{equation}
In this way eqn.~\eqref{vlasov1} can be solved with the method of characteristics.
Given the fields $B$ and $E$ we can follow the characteristics by solving eqn. \eqref{vlasov1:char}
with a standard ODE integrator.
We can introduce a second density $g(x,v,t)$ which solves the same Vlasov equation as $f$
\begin{equation}
\label{vlasov2}
  \partial_t g(x,v,t) + v \cdot \partial_x g(x,v,t) + \frac{q}{m}E(x,t) \cdot \partial_v g(x,v,t) =0
 \end{equation}
and call it the sampling density, prior or the law of $(\mathbf{X},\mathbf{V})$.
The initial sampling distribution $g(\cdot,\cdot,t=0)$ becomes a probability density
by imposing a normalization over the phase space $\Omega$ by $\int_{\Omega} g(x,v,t=0) \dxdv =1$ and $g(x,v,t=0)\geq 0$ for all $(x,v) \in \Omega$ .
Since $g$ follows the same Vlasov equation \eqref{vlasov2}
as $f$, see eqn.~\eqref{vlasov1}, it is constant along the same characteristics \eqref{vlasov:characteristics}.\\
The Vlasov equation $\eqref{vlasov2}$ conserves positivity and volume, therefore, $g$ stays a probability density for all $t\geq 0$, which is discussed and verified in the next section.\\
\subsection{Phase space conservation}
In order to verify that $g(x,v,t)$ is the probability density of the characteristics $(\mathbf{X}(t),\mathbf{V}(t))$ we rewrite the characteristics as a mapping $\varphi_t$.
Since $f$ is constant along the characteristics, we can implicitly define a diffeomorphism $\varphi_t: (x_0,v_0) \mapsto (x,v) $ for every $t \geq 0$ such that
\begin{align}
 \label{vlasov:characteristics:flux1}
 f(x,v,t)=f(\varphi_t(x_0,v_0),t) = f(x_0,v_0, 0).
\end{align}
The same property then also holds for $g$, namely $g(\varphi_t(x_0,v_0),t) = g(x_0,v_0, 0)$. We seek a change in variables $(x,v):=\varphi_t(x_0,v_0)$, as we are interested
in what happens with $f$ and $g$ at later times. For this denote the Jacobi determinant of $\varphi_t$ as $J_{\varphi_t}$. 
In general, after a transformation has been applied onto a random deviate it's probability density has to be scaled with the according Jacobi determinant, as we will recall in the next step.
For any phase-space volume $V\subset{\Omega}$ equation~\eqref{vlasov:characteristics:flux2} then holds under the change of variables; also for $f$.
\begin{equation}
 \label{vlasov:characteristics:flux2}
 \begin{split}
 \iint_{\varphi(V)} g(x,v,t) \dxdv
 &=  \iint_{V}  g\left(\varphi_t(x_0,v_0) ,t\right)  J_{\varphi_t}(x_0,v_0) ~\mathrm{d}x_0 \mathrm{d}v_0\\
 &= \iint_{V}  g\left(x_0,v_0 ,0\right)  J_{\varphi_t}(x_0,v_0) ~\mathrm{d}x_0 \mathrm{d}v_0
 \end{split}
\end{equation}
This means that at time $t$, $(x_0,v_0)\mapsto g\left(x_0,v_0,t=0 \right) \cdot J_{\varphi_t}(x_0,v_0)$ is the probability density
for the random deviate $(\mathbf{X}(t),\mathbf{V}(t))=\varphi(X_0,V_0)$ and the Jacobian has to be taken into account. 
For the Vlasov equation the Jacobi determinant is one, $J_{\varphi_t}(x,v)=1$. Hence the characteristics transport
the actual value of the probability density at every time $t$. This also holds true for a symmetric integrator, e.g.~one time step of the
symplectic Euler scheme given in equation~\eqref{vlasov:characteristics:flux:seuler}.
\begin{equation}
\label{vlasov:characteristics:flux:seuler} 
\begin{split}
\varphi_t(x,v) =  \left(x + tv,  ~v + t \frac{q}{m} E(x+tv,0) \right),\\
\nabla \varphi_t(x,v)=
\begin{pmatrix}
 1 & t \\
 t \frac{q}{m} \partial_x E(x+tv,0) &  1+ t^2 \frac{q}{m} \partial_x E(x+tv,0) \\
\end{pmatrix}
\end{split}
\end{equation}
We then see that the semi-discrete flow also has the right Jacobi determinant:
\begin{equation}
\det(\nabla \varphi_t) = 1+ t^2 \frac{q}{m} \partial_x E(x+tv,0) - t^2 \frac{q}{m} \partial_x E(x+tv,0)  = 1.
\end{equation}
Yet when we consider the standard explicit Euler scheme and its Jacobi determinant
given in eqn.~\eqref{vlasov:characteristics:flux:euler} the determinant of the flow is not one. 
\begin{equation}
\label{vlasov:characteristics:flux:euler} 
\begin{split}
\varphi_t(x,v) &=  \left(x + tv,  ~v + t \frac{q}{m} E(x,0) \right),
\quad
\nabla \varphi_t=
\begin{pmatrix}
 1 & t \\
 t \frac{q}{m} \partial_x E(x,0) &  1\\
\end{pmatrix}
\end{split}
\end{equation}
\begin{equation}
\label{explicit:euler:jacobideterminant}
\det(\nabla \varphi_t) = 1 - t^2 \frac{q}{m} \partial_x E(x,0) \neq1
\end{equation}
Therefore the likelihood $g$ has to be rescaled accordingly such that it continuously represents the 
distribution of the random deviate $(\mathbf{X}(t),\mathbf{V}(t))$. Technically, $f$ should still stay constant because the we use the method of characteristics, which leads ultimately to an inconsistency.\\
By symmetric composition it is possible to extend the dissipative explicit Euler scheme, given by $\varphi_t$ with it's left adjoint defined over the time inverted inverse $\varphi^*_t = \varphi_{-t}^{-1}$.
\begin{equation}
 \det(\nabla \varphi_t^*) = \det(\nabla \varphi_{-t}^{-1})
= \frac{1}{ \det(\nabla \varphi_{-t}) } = \frac{1}{ 1 - (-t)^2 \frac{q}{m} \partial_x E(x,0) } = \frac{1}{ \det(\nabla \varphi_{t}) } 
\end{equation}
Therefore, the composition of the explicit Euler with it's adjoint should provide us a phase space conserving method. Instead of one full time step we move only a half step
and obtain a second order method. In order to work with the discrete mappings, we introduce the discrete time grid by $t_n= n\, \Delta t$. For the Vlasov--Poisson system the
electric field is obtained by the position of the particles, so it is important to note that the electric field $E(x,t_n)$ is determined from the particles $x_n$ at the $n^{\mathrm{th}}$ time step.
The discrete explicit Euler and it's adjoint read then:
\begin{align}
\label{discrete:explicit:euler}
 \varphi_{ \frac{\Delta t}{2} }&
\begin{cases}
  x_{n+1} =  x_{n+ \nicefrac{1}{2}} + \frac{\Delta t}{2} v_{n+ \nicefrac{1}{2}}\\
  v_{n+1} =  v_{n+ \nicefrac{1}{2}} + \frac{\Delta t}{2} \frac{q}{m}  E(x_{n+\nicefrac{1}{2}},t_{n+\nicefrac{1}{2}} )\\
\end{cases}\\
\label{discrete:implicit:euler}
\varphi^*_{\frac{\Delta t}{2}}= \varphi^{-1}_{-\frac{\Delta t}{2}}&
=\begin{cases}
  x_{n+ \nicefrac{1}{2}} =   x_n + \frac{\Delta t}{2} v_{n+\nicefrac{1}{2}}\\
  v_{n+ \nicefrac{1}{2}} =  v_n + \frac{\Delta t}{2} \frac{q}{m}  E( x_n + \frac{\Delta t}{2} v_{n+\nicefrac{1}{2}},t_{n+\nicefrac{1}{2}} )\\
 \end{cases}
  \end{align}
The two ways of combining \eqref{discrete:explicit:euler} and \eqref{discrete:implicit:euler} are the Crank--Nicolson,
\begin{equation}
\label{discrete:cranknicolson}
  \varphi^*_{\frac{\Delta t}{2}} \circ \varphi_{ \frac{\Delta t}{2} }
 \begin{cases}
  x_{n+1} =  x_n +  \frac{\Delta t}{2} \, \frac{v_n + v_{n+1} }{2}\\
  v_{n+1} =   v_n + \frac{\Delta t}{2} \, \frac{q}{m} \left[ E( x_n ,t_{n} ) + E( x_{n+1} ,t_{n+1} )    \right] \\
 \end{cases}
\end{equation}
and the implicit midpoint method:
\begin{equation}
\label{discrete:implicitmidpoint}
  \varphi_{ \frac{\Delta t}{2} } \circ \varphi^*_{\frac{\Delta t}{2}} 
 \begin{cases}
  x_{n+ \nicefrac{1}{2}} =   x_n + \frac{\Delta t}{2} v_{n+\nicefrac{1}{2}}\\
  v_{n+ \nicefrac{1}{2}} =  v_n + \frac{\Delta t}{2} \frac{q}{m}  E( x_n + \frac{\Delta t}{2} v_{n+\nicefrac{1}{2}},t_{n+\nicefrac{1}{2}} )\\
  x_{n+1} =  x_n +  \Delta t \, v_{n+ \nicefrac{1}{2}}\\
  v_{n+1} =   v_n + \Delta t \, \frac{q}{m}  E( x_n + \frac{\Delta t}{2} v_{n+\nicefrac{1}{2}},t_{n+\nicefrac{1}{2}} )\\
 \end{cases}
\end{equation}
Both methods are fully implicit, but the crucial difference is that the implicit midpoint method approximates the electric field in the middle of the time step, whereas
the Crank--Nicolson averages the field at the beginning and the end of each time step. While both methods are considered geometric integrators with excellent long term stability
only one of them conserves the phase space volume exactly. We have already determined in eqn.~\eqref{vlasov:characteristics:flux:euler} that the explicit Euler did not conserve phase space and is dissipative,
hence it's adjoint has to be investigated. The map underlying the implicit Euler \eqref{discrete:implicit:euler} reads
\begin{equation}
 \varphi^*_t(x,v) = \left(\varphi^*_{t,x}(x,v)  ,\varphi^*_{t,v}(x,v) \right) = \left(x + t \varphi^*_{t,x}(x,v),  ~v + t \frac{q}{m} E(\varphi^*_{t,x}(x,v) ,t) \right)  
 \end{equation}
Although the map is given implicitly, the Jacobi matrix can be calculated by straight forward derivation and yields also an implicit expression:
 \begin{multline}
 \nabla \varphi^*_t(x,v)=
 \begin{pmatrix}
  \partial_x \varphi^*_{t,x} & \partial_v \varphi^*_{t,x}\\
  \partial_x \varphi^*_{t,v} & \partial_v \varphi^*_{t,v}
 \end{pmatrix}\\
=
 \begin{pmatrix}
  1 + t \partial_x \varphi^*_{t,v}(x,v) & t \partial_v \varphi^*_{t,v}(x,v) \\
  t \frac{q}{m} \partial_x E( \varphi^*_{t,x}(x,v), t )  \partial_x \varphi^*_{t,x}(x,v) 
  &  1 +  t \frac{q}{m} \partial_x E( \varphi^*_{t,x}(x,v), t ) \partial_v \varphi^*_{t,x}(x,v) \\
 \end{pmatrix} 
\end{multline}
By the suitable insertion of the implicitly defined derivatives the Jacobi determinant reduces to
\begin{multline}
\det(\nabla \varphi^*_t) = 
    \underbrace{ \left[ 1 + t \partial_x \varphi^*_{t,v}(x,v) \right]}_{= \partial_x \varphi^*_{t,x}(x,v)}
    \left[ 1 +  t \frac{q}{m} \partial_x E( \varphi^*_{t,x}(x,v), t )  \underbrace{ \partial_v \varphi^*_{t,x}(x,v)}_{=t \partial_v \varphi^*_{t,v}(x,v)} \right]\\
  -t \frac{q}{m} \partial_x E( \varphi^*_{t,x}(x,v), t )  \partial_x \varphi^*_{t,x}(x,v) t \partial_v \varphi^*_{t,v}(x,v)= \partial_x \varphi^*_{t,x} (x,v).
\end{multline}
Unfortunately $\partial_x \varphi^*_{t,x} (x,v)$ is not known to us, but re-substituting expressions from the Jacobi matrix yields a recurrence relation
\begin{equation}
 \partial_x \varphi^*_{t,x}(x,v)=  1 + t \partial_x \varphi^*_{t,v}(x,v)  = 1+ t^2 \frac{q}{m} \partial_x E( \varphi^*_{t,x}(x,v), t )  \partial_x \varphi^*_{t,x}(x,v),
\end{equation}
which is easily resolved:
\begin{multline} 
 \Rightarrow
  1 =  \partial_x \varphi^*_{t,x}(x,v) - \partial_x \varphi^*_{t,x}(x,v) t^2 \frac{q}{m} \partial_x E( \varphi^*_{t,x}(x,v), t )    \\
  \Rightarrow 
  \partial_x \varphi^*_{t,x}(x,v)  = \frac{1}{1-\partial_x \varphi^*_{t,x}(x,v) t^2 \frac{q}{m} \partial_x E( \varphi^*_{t,x}(x,v), t )    }\\
  \Rightarrow
  \det(\nabla \varphi^*_t) = \frac{1}{1 - t^2 \frac{q}{m} \partial_x E( \varphi^*_{t,x}(x,v), t )}
  \label{implicit:euler:jacobideterminant}
\end{multline}
We realize that \eqref{implicit:euler:jacobideterminant} constitutes the inverse \eqref{explicit:euler:jacobideterminant} if and only if the 
implicit Euler is applied first, as in \eqref{discrete:implicitmidpoint}. Note that in this case $\varphi^*_{\nicefrac{\Delta t}{2},x}(x_n,v_n) = x_{n+\nicefrac{1}{2}}$, which means 
that the Jacobi determinant of the implicit midpoint scheme is one:
\begin{multline}
 \det\left[ \nabla \Big( \varphi_{\nicefrac{\Delta t}{2}} \circ \varphi^*_{ \nicefrac{\Delta t}{2}} (x_n,v_n)  \Big) \right]
 =\det\left[ \nabla \varphi_{\nicefrac{\Delta t}{2}} \big( \underbrace{\varphi^*_{ \nicefrac{\Delta t}{2}} (x_n,v_n)}_{ =( x_{n+\nicefrac{1}{2}},  v_{n+\nicefrac{1}{2}}) } \big) 
  \nabla \varphi^*_{ \nicefrac{\Delta t}{2}} (x_n,v_n) \right]\\
 =   \frac{1 - \big(\frac{\Delta t}{2} \big)^2 \frac{q}{m} \partial_x E\big(  x_{n+\nicefrac{1}{2}}, t_{n+\nicefrac{1}{2}} \big) 
  }{   1 - \big(\frac{\Delta t}{2} \big)^2 \frac{q}{m} \partial_x E( \varphi^*_{\nicefrac{\Delta t}{2},x}(x_n,v_n), t_{n+\nicefrac{1}{2}} )   } 
  =   \frac{1 - \big(\frac{\Delta t}{2} \big)^2 \frac{q}{m} \partial_x E\big(  x_{n+\nicefrac{1}{2}}, t_{n+\nicefrac{1}{2}} \big) 
  }{   1 - \big(\frac{\Delta t}{2} \big)^2 \frac{q}{m} \partial_x E(x_{n+\nicefrac{1}{2}} , t_{n+\nicefrac{1}{2}} )   } =1
\end{multline}
Contrary, for the Crank--Nicolson this holds not true. The Jacobi determinants cancel out in between half time steps, such that any series of time steps begins and ends
with a slightly dissipative half step. Hence one can call the Crank-Nicolson \textit{adjoint} phase space conserving.\\
It is important to note that the Vlasov--Poisson system is a Hamiltonian system in which our phase space coordinates $(x,v)$ coincide
with the Hamiltonian coordinates $(q,p)$. Without magnetic field the system for a single partile can be written as
\begin{equation}
 (\dot{p},\dot{q})= J^{-1} \nabla_{(p,q)} H(p,q), \quad
 J= \begin{pmatrix}
  & -I\\
  I&\\
 \end{pmatrix},
\end{equation}
with $H(p,q) = \frac{p^2}{2} + \Phi(q)$. 
For different systems, we will obtain a different matrix $J$ and the coordinates $(p,q)$ cannot be identified as $(x,v)$ much longer.
An integrator is called symplectic if the mapping induced by $\varphi_t$
is symplectic with respect to $J$, which is checked by
\begin{equation}
\nabla \varphi_t(p,q)^t J \nabla \varphi_t(p,q) = J.
\end{equation}
See Hairer's lecture notes for a short introduction to Hamiltonian systems~\cite{hairer_tum2010}. Such symplectic integrators
always conserve phase-space volume and can also conserve quantities like energy but not every phase space volume-preserving integrator is symplectic, see also~\cite{mclachlan2006geometric}. 
But conservation of phase space is such an important property that schemes like the Boris method perform so well although they cannot be symplectic for any system~\cite{hongqin2013}.
Many of these integrators along with detailed theory for plasma physics can already be found in~\cite{kraus2013variational}. For the Vlasov--Poisson system the commonly known schemes are symplectic Runge Kutta methods up to fourth order~\cite{forest1989fourth}, where
second order scheme corresponds to the well-known leap frog, and the first order is the symplectic Euler.
But so far it is unclear whether symplecticity provides advantages concerning the conservation of low-discrepancy, which is why we restrict ourselves to phase space conservation.

\subsection{Monte Carlo integration with particles}
So far we did not address how to solve the actual Poisson equation and acquire the electric fields.
We will now slightly deviate in notation from the standard Particle-In-Cell (PIC) method~\cite{birdsall2004plasma}. The introduction of the probability density function $g$ allows us to define the characteristics $\mathbf{X}(t)$ and $\mathbf{V}(t)$ as random variables for each time $t$, such that the 
trajectories in time form a stochastic process~\cite{oksendal2003stochastic}
describing the solution to eqn.~\eqref{vlasov1}.\\
Before solving the Poisson equation one can explain the stochastic setting by estimating the kinetic energy $\mathcal{H}_T$ which is a moment of the solution $f$. The characteristics are random variables with joint probability density $g$. Suppose $g$ supports $f$ which means $ \mathrm{supp}(f) \subset \mathrm{supp}(g) $. Then by inserting $g$ integrals over $f$ can be linked to expected values over $X$ and $V$.
\begin{equation}
\label{estimator:ekin1}
\begin{split}
\mathcal{H}_T(t) &= \frac{1}{2} \int_{\Omega} v^2 \, f(x,v,t) \dxdv \\
& = \frac{1}{2} \int_{\Omega} v^2 \, \frac{f(x,v,t)}{g(x,v,t)}~g(x,v,t) \dxdv \\ 
 &=\frac{1}{2} \,\mathbb{E}\left[ \mathbf{V}(t)^2 \frac{f(\mathbf{X}(t),\mathbf{V}(t),t)}{g(\mathbf{X}(t),\mathbf{V}(t),t)} \right] \\
 \end{split}
\end{equation}
The values of $g$ and $f$ over time are constant along the characteristics (see eqn.~\eqref{vlasov:characteristics})
such that eqn.~\eqref{estimator:ekin1} simplifies to eqn.~\eqref{estimator:ekin2}.
\begin{equation}
\label{estimator:ekin2}
\mathcal{H}_T(t) =
\frac{1}{2} \,\mathbb{E}\left[ \mathbf{V}(t)^2 \frac{f(\mathbf{X}(t),\mathbf{V}(t),t)}{g(\mathbf{X}(t),\mathbf{V}(t),t)} \right] 
=\frac{1}{2}\, \mathbb{E}\left[ \mathbf{V}(t)^2 \frac{f(X(0),V(0),0)}{g(X(0),V(0),0)} \right] \\
\end{equation}
In order to get an estimate of the expectation in eqn.~\eqref{estimator:ekin2} one has to use Monte Carlo integration. We define $N_p$ independently and identically distributed (i.i.d.) samples 
$\left(\mathbf{x}^0_k, \mathbf{v}^0_k \right)_{k=1,\dots,N_p}$ of the random deviates $(X(0),V(0))$ using the knowledge of the probability density $g(x,v,t=0)$. These samples are called markers or particles.
The samples $\left(\mathbf{x}^t_k, \mathbf{v}^t_k \right)_{k=1,\dots,N_p}$ can be advanced over time using a suitable phase space conserving time integrator. Then at any point in time they are distributed according to $g(\cdot,\cdot, t)$ as a solution to the Vlasov equation.
Note that the plasma likelihood
\begin{equation}
\mathbf{f}^t_k := f( \mathbf{x}_k^t, \mathbf{v}_k^t,t) = f( \mathbf{x}_k^0, \mathbf{v}_k^0,0)      =\mathbf{f}^0_k 
\end{equation}
and the sampling likelihood
\begin{equation}
\mathbf{g}^t_k := g( \mathbf{x}_k^t, \mathbf{v}_k^t,t) = g( \mathbf{x}_k^0, \mathbf{v}_k^0,0) \mathbf{g}^0_k 
\end{equation}
stay constant over time. In the common notation of collisionless PIC schemes
the ratio between those two likelihoods is referred to as the time-independent particle weight $\mathbf{w}_k^0$.
\begin{equation}
 \mathbf{w}_k= \mathbf{w}_k^t = \frac{\mathbf{f}^t_k }{\mathbf{g}^t_k } = \frac{\mathbf{f}^0_k }{\mathbf{g}^0_k } =  \mathbf{w}_k^0
\end{equation}
This allows us to estimate the kinetic energy using the samples in eqn.~\eqref{estimator:ekin3}.
\begin{equation}
\label{estimator:ekin3}
\begin{split}
\mathcal{H}_T(t) &=\frac{1}{2}\, \mathbb{E}\left[ \mathbf{V}(t)^2 \frac{f(X(0),V(0),0)}{g(X(0),V(0),0)} \right]\\
&\approx \frac{1}{2}\,\frac{1}{N_p} \sum_{k=1}^{N_p} \big(\mathbf{v}_k^t\big)^2 \frac{\mathbf{f}^0_k }{\mathbf{g}^0_k } 
= \frac{1}{2}\, \frac{1}{N_p} \sum_{k=1}^{N_p} \big(\mathbf{v}_k^t\big)^2 \mathbf{w}_k.
\end{split}
\end{equation}
Note that sometimes the factor $\frac{1}{N_p}$ is included into the particle weight $\mathbf{w}_k$.
\subsection{Weak Poisson solve with particles}
It only remains to solve the Poisson equation~\eqref{poisson1} using the samples $\left(\mathbf{x}^0_k, \mathbf{v}^0_k \right)_{k=1,\dots,N_p}$. This is commonly done in weak form given a test function $\text{ for } \varphi \in H_1([0,L])$, see also \cite{kraus2016gempic}.
\begin{equation}
 - \int_0^L \Phi(x,t)  \varphi(x) \dx  
 = \int_0^L q\Big(\int_{-\infty}^{\infty} f(x,v,t)   \mathrm{d}v -1 \Big) \varphi(x) \dx
\end{equation}
The only unknown is the right hand side, 
\begin{equation}
q \int_0^L \int_{-\infty}^{\infty} f(x,v,t)   \varphi(x) \,\mathrm{d}v \mathrm{d} x
\end{equation}
which depends on the density $f$. But in the same manner as we calculated the 
kinetic energy we can use Monte Carlo estimator to calculate the right hand side.
\begin{equation}
\begin{split}
q \int_0^L \int_{-\infty}^{\infty} f(x,v,t)   \varphi(x) \,\mathrm{d}v \mathrm{d} x
&=  q\,\mathbb{E} \left[  \varphi\big(\mathbf{X}(t)\big) \frac{f(\mathbf{X}(t),\mathbf{V}(t),t)}{g(\mathbf{X}(t),\mathbf{V}(t),t)} \right]\\
&\approx
 q\,\frac{1}{N_p} \sum_{k=1}^{N_p} \big(\mathbf{v}_k^t\big)^2 
\underbrace{\frac{\mathbf{f}^0_k }{\mathbf{g}^0_k } }_{=\mathbf{w}^0_k}
\end{split}
\end{equation}
In order to obtain the common Particle-in-Cell formulation a spline basis as in eqn.~\eqref{linearsplinebasis} can be used for the test function $\varphi$ and the solution $\Phi$.
The electric field is then easily obtained $E(x,t)= - \nabla \Phi(x,t)$ at any time such that it can be used in the phase-space conserving ODE integrators for advancing the particles.

\section{Monte Carlo and Quasi Monte Carlo}
\begin{figure}[H]  
\centering
 \includegraphics[width=0.8\textwidth]{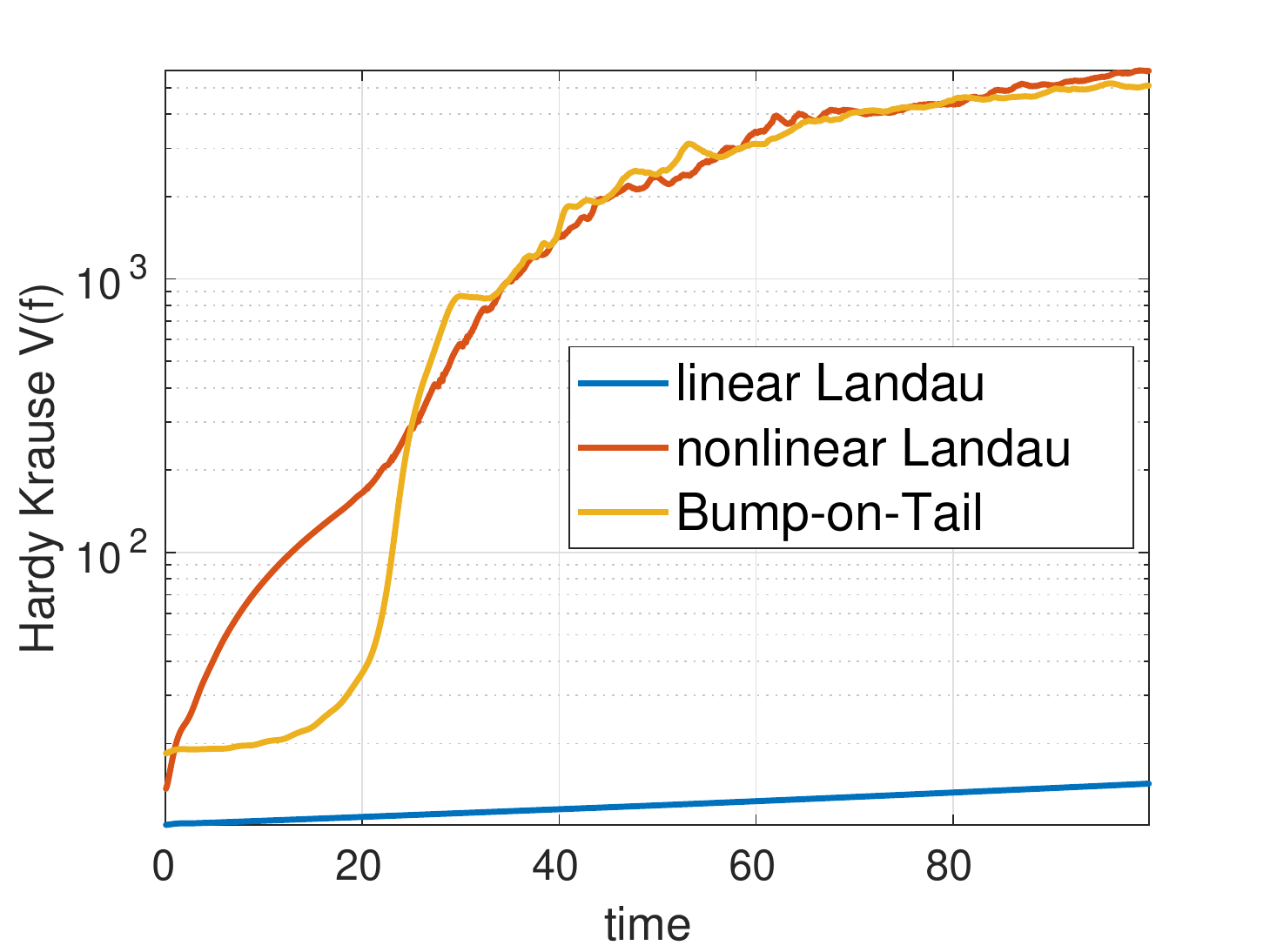}
 \caption{Total variation of the distribution function $f$ in the sense of Hardy and Krause for different test cases of the Vlasov Poisson. Using 
           pseudo spectral solver with a third order symplectic Runge Kutta scheme at $\Delta t = 0.05$ and high resolution $N_x=512, ~N_v=1024$ the total variation
           is calculated at each time step according to eqn.~\eqref{tvhk2d}. For linear Landau damping, there are only minor changes in the distribution function since only the recurrence phenomenon is causing a
           long time disturbance. Both nonlinear Landau damping and the Bump-on-tail instability create many small scale features in the distribution function, which leads to the harsh increase in variation.
           The vortex in the Bump-on-tail instability, also known as the BGK mode, take some time to fully develop whereas the nonlinear Landau damping ``folds''
           the phase space much quicker, which explains the difference in the initial development of the variation.}
 \label{fig:tv_hk_overview}
\end{figure}

\begin{figure}[H]
\centering
\begin{tabular}{c | c}
 \hline
 \multicolumn{2}{c}{Pseudo random}\\
 \includegraphics[width=0.48\textwidth]{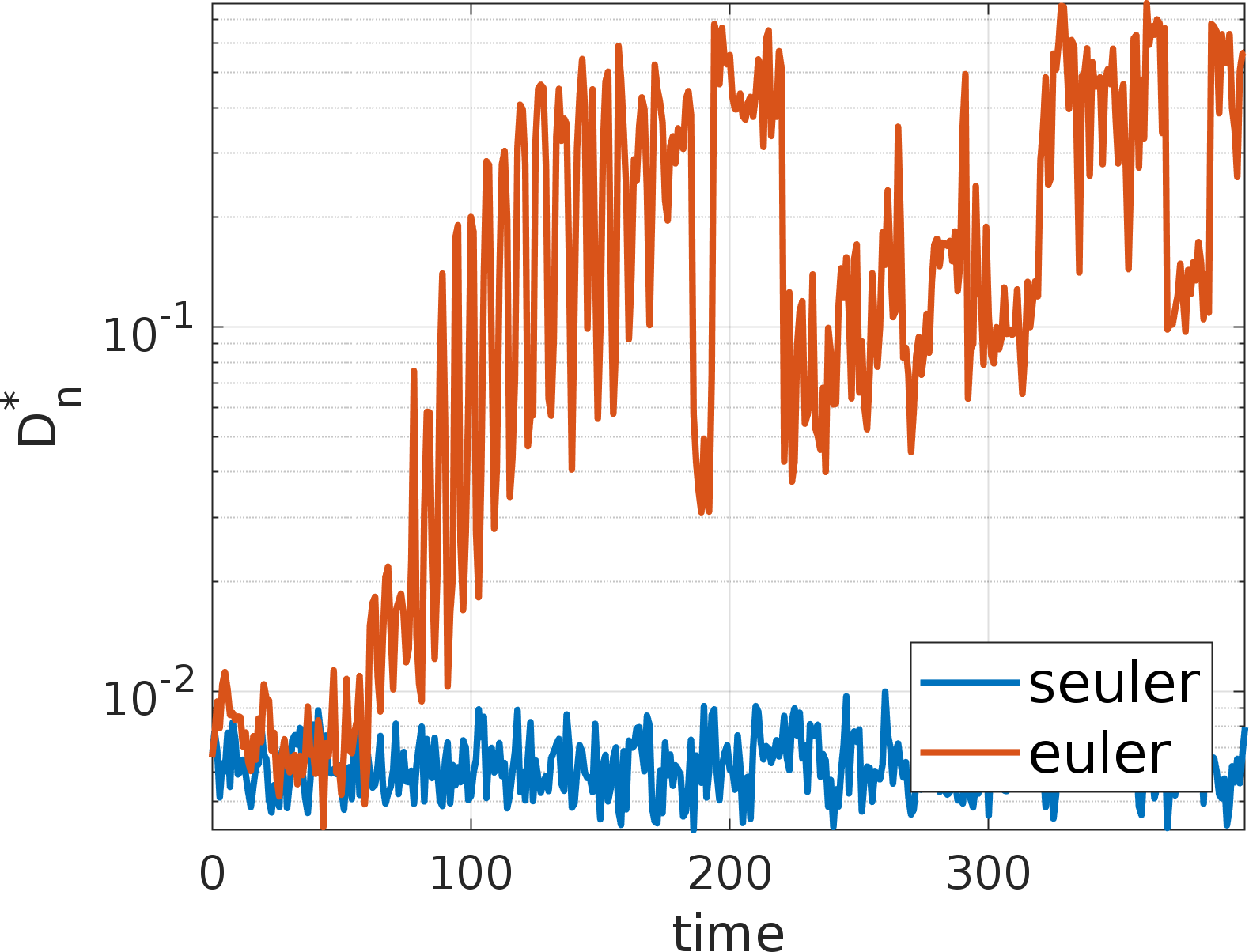}
 &
 \includegraphics[width=0.48\textwidth]{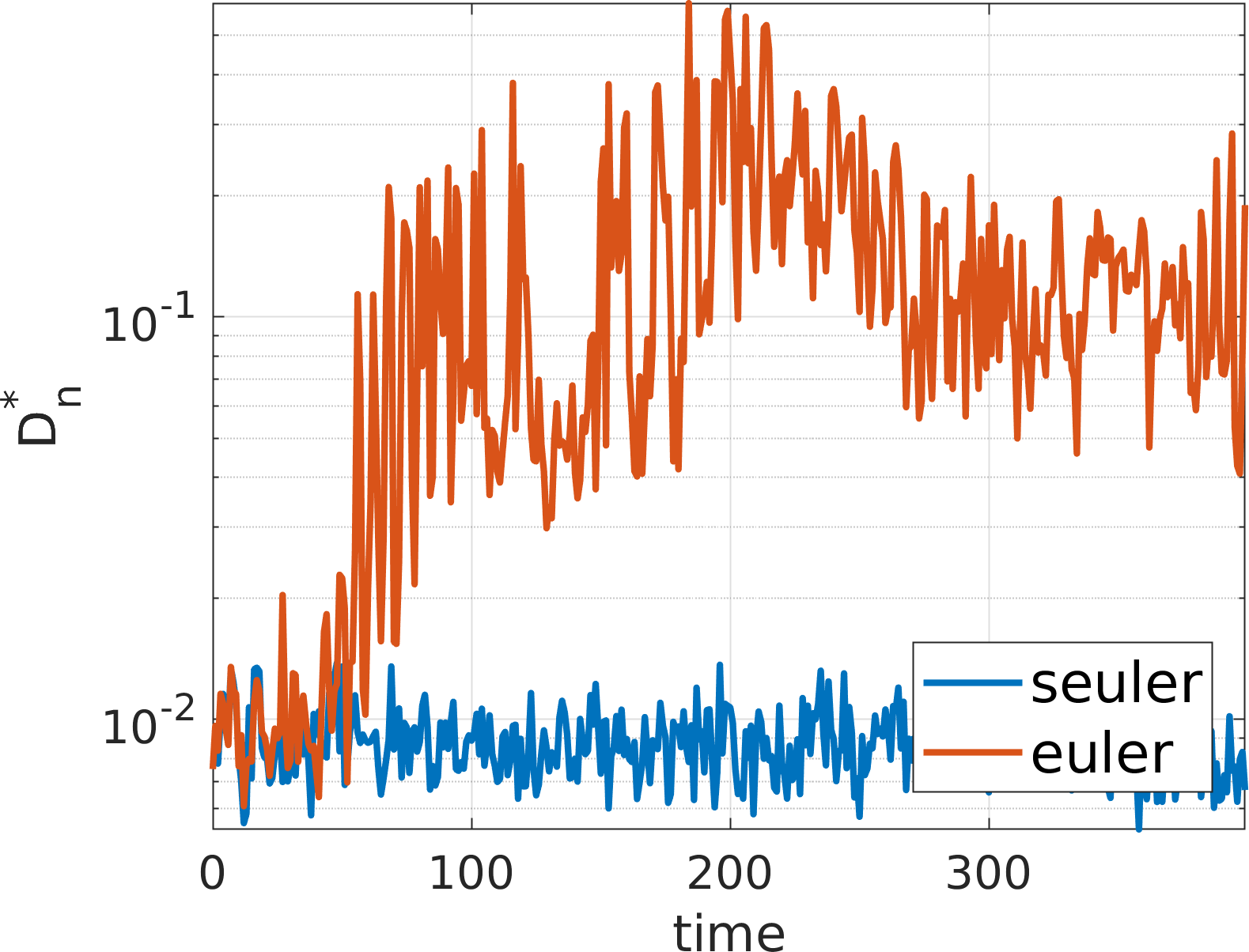}\\
 \hline
  \multicolumn{2}{c}{Quasi Monte Carlo (Sobol)}\\
 \includegraphics[width=0.48\textwidth]{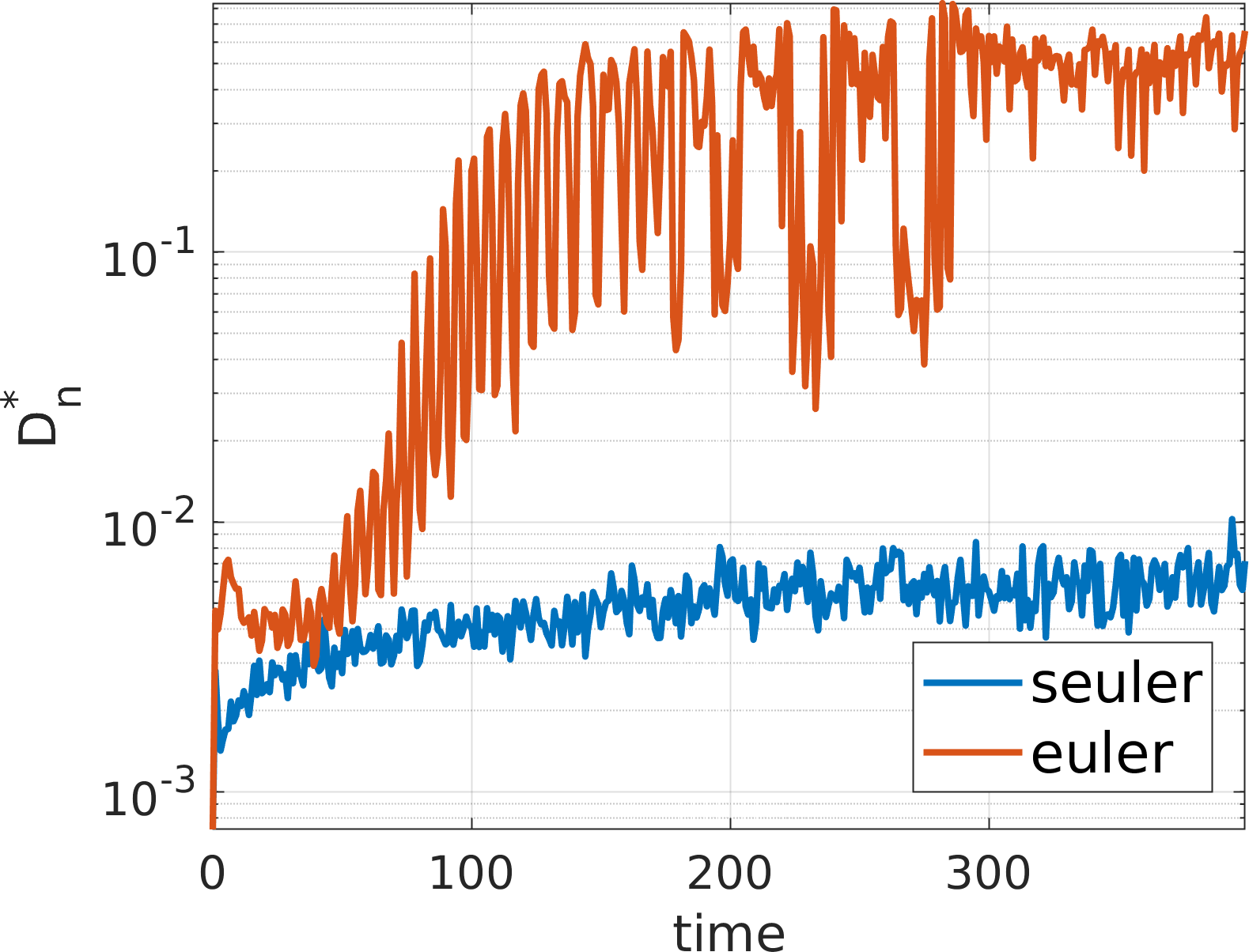}&
 \includegraphics[width=0.48\textwidth]{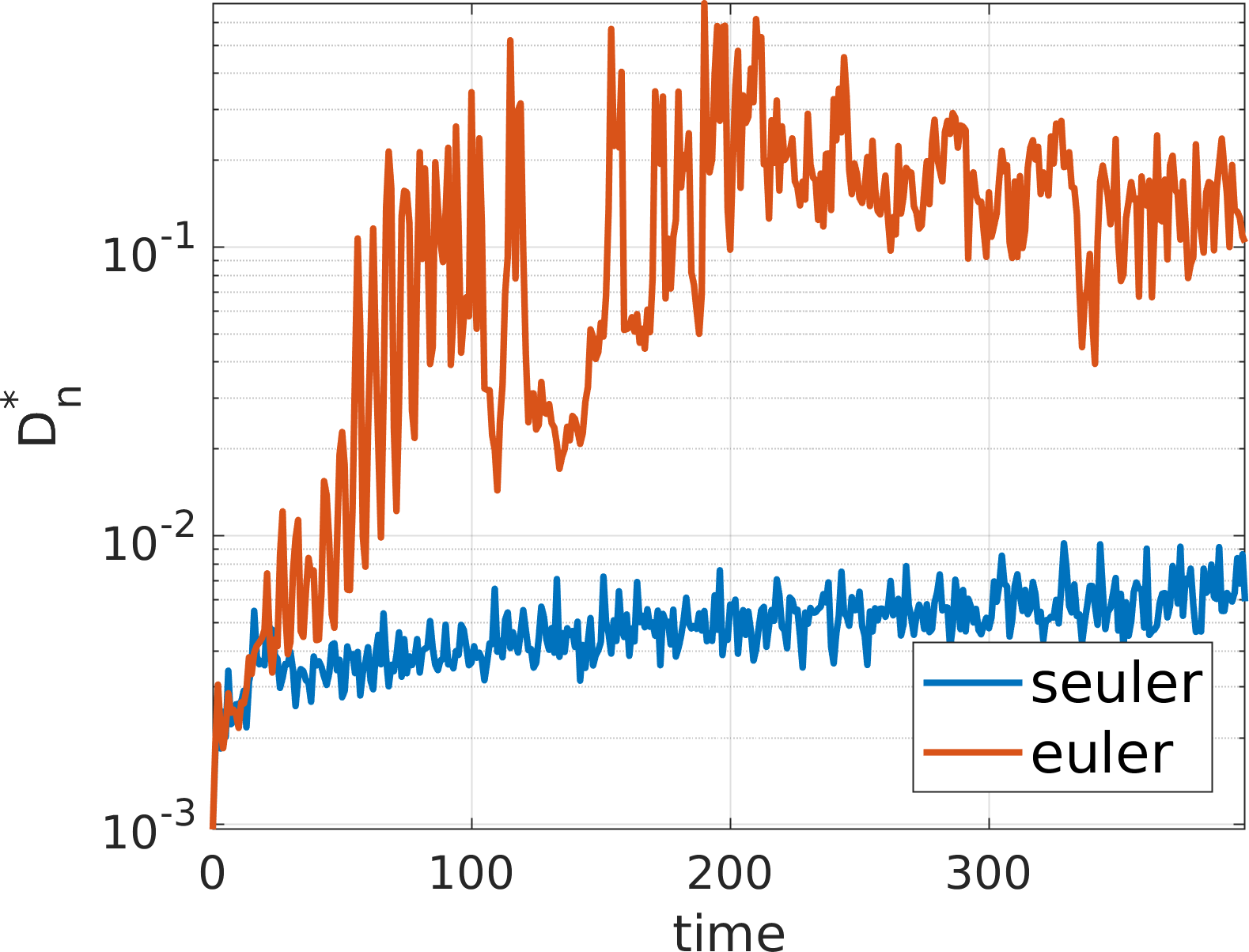}\\
  \hline
 (a) nonlinear Landau damping  & (b) Bump-on-tail instability\\
 \end{tabular}
\caption{The star discrepancy in a long term PIC simulation of nonlinear Landau damping (a) and a Bump-on-tail instability (b) with a total of $N_p=2\cdot10^6$ particles sampled uniformly from the Sobol sequence. The explicit euler ``euler'' 
          is dissipative which worsens the discrepancy compared to the phase space conserving symplectic euler ``seuler''.}
\label{fig:pic:star_discrepancy}
\end{figure}

As commonly known the Monte Carlo estimator converges with $\mathcal{O}(N^{-0.5})$ when using $N$ pseudo random samples. Samples obtained from low discrepancy sequences,
also known as Quasi Monte Carlo (QMC), can obtain convergence up to $\mathcal{O}(N^{-1} \log(N)^{d-1}  )$ under certain smoothness conditions onto the integrand~\cite{morokoff1994quasi,basu2016transformations}. Also, the convergence rate is not restricted to the Lesbegue measure~\cite{aistleitner2016tusnady}.
Yet the measure of error for integration with these low discrepancy sequences is the Hardy-Krause variation, see~\cite{hickernell2005} for an overview and also \cite{drmota2006sequences,hlawka1971definition}.
High order scrambling by Dick, see~\cite{dick2011higher,dick2013high}, leads to convergence rates up to $\frac{7}{2}$ but requires even smoother integrands.

\begin{equation}
\label{eq:Koksma-Hlawka}
 \Big| \int f(x) \mathrm{d}x - \frac{1}{n} \sum_{p=1}^n f(x_p)   \Big| \leq V(f) D_n^* 
\end{equation}
where $V(f)$ denotes the Hardy Krause variation of $f$ and $D_n^*$ the star discrepancy of the point set $x_n$. 
The total variation in the sense of Hardy Krause depends on the scale of the smallest features of $f$ and can be calculated for continuously differentiable functions, see \cite{morokoff1994quasi}.
Almost any randomly chosen sequence will be bounded as
\begin{equation}
  D_n = \mathcal{O} \Bigg( \sqrt{ \frac{ \log( \log(n ))}{n}} \Bigg) \text{ with } D_n^* \leq  D_n \leq 2^d D_n^*,
\end{equation}
hence the observed $\mathcal{O}\sqrt{n}$ convergence. On the other hand the Halton and Sobol quasi random sequences yield a asymptotically smaller star discrepancy:
\begin{equation}
 D_n^* = C_d \frac{\log(n)^d}{n} + \mathcal{O}\Bigg( \frac{ \log(n )^{d-1}}{n}  \Bigg) ,  \quad C_d \ge 0.
\end{equation}
For an overview over different discrepancies and sequences~\cite{morokoff1994quasi} is recommended. For the Sobol sequence the best upper bound
for $C_d$ is $5.28$. But in~\cite{thiemard2001optimal} it is shown that the current theory on those bounds is insufficient for any practical use. Note that the total variation is defined as
\begin{equation}
 V(f) = \int | \nabla f(x)  | \mathrm{d} x,
\end{equation}
but in the world of low discrepancy sequences one mostly uses the definition in the sense of Hardy and Krause, which is the sum of the $L^1$ norm of all first order partial derivatives, except that
the identically mixed derivatives are only accounted once. For our two dimensional phase space this reads
\begin{equation}
\label{tvhk2d}
  V(f) =  \iint | \partial_{x} f(x,v,t) | \, \mathrm{d}v\mathrm{d}x 
  \,+\,  \iint | \partial_{v} f(x,v,t) | \, \mathrm{d}v\mathrm{d}x
  \,+\,  \iint | \partial_{x} \partial_{v} f(x,v,t) | \, \mathrm{d}v\mathrm{d}x. 
\end{equation}
Naturally, the question arises how the total variation changes over time in the Vlasov--Poisson system. Since it is quite hard to compute~\eqref{tvhk2d} from a marker distribution we used a
pseudo spectral solver. Here fig.~\ref{fig:tv_hk_overview} clearly shows increases in the variation of several orders of magnitudes due to the development of small scales. This means
according to the Koksma-Hlawka inequality~\ref{eq:Koksma-Hlawka}, that even if the discrepancy of the markers stays constant much more markers are needed in the nonlinear phase.
Therefore also the star discrepancy of the markers in a PIC simulation has to be investigated more closely.
There are various ways of estimating the (star) discrepancy of a given point set~\cite{doerr2014calculation}, where we chose the method provided by~\cite{thiemard2001optimal}.
The original implementation provided by~\cite{thiemard2001optimal} is restricted to computing the discrepancy in a quadratic box. Hence we periodically estimate the 
star discrepancy of all markers in the phase space box $ (x,v) \in [0,2]\times[-1,1] $. Originally the markers are sampled uniformly in phase space, such that we expect this
uniformity to be preserved over time in any subinterval. This way, unfortunately, the sampling has to be cut off in the velocity domain such that we chose $|v|\leq 8$ initially.
The testing box is chosen smaller to avoid the influence of this boundary since it is not guaranteed that the support of the sampling distribution $g$ stays constant.
We already know that the standard Euler is dissipative, so we expect it to influence the star discrepancy. As can be seen fig.~\ref{fig:pic:star_discrepancy}
the symplectic Euler suffers only from a minor change in the discrepancy whereas the standard Euler suffers from an increase of orders of magnitude.
In view of Koksma-Hlawka inequality~\eqref{eq:Koksma-Hlawka} these results strongly recommend the use of uniformity preserving methods, such as the symplectic Euler because otherwise much more markers are required.

\subsection{Inverse Transform Sampling}
For a given phase space density $f(x,v)$ we can define the sampling density $g$ as
\begin{equation}
\label{normalizef}
 g_{\mathbf{X},\mathbf{V}}(x,v)= \frac{|f(x,v)|}{ \int_{ v_{\min}}^{v_{\max}} \int_{x_{\min}}^{x_{\max}} |f(x,v)|~  \mathrm{d}x\mathrm{d}v }
\end{equation}
There are various ways of sampling from an arbitrary probability density $g$. For pseudo random numbers popular choices are Markov Chain Monte Carlo, Gibbs sampling or even 
the inefficient rejection sampling~\cite{caflisch98,liu2008monte,owen2013,christian2007monte}. 
Although there exist MCMC algorithms~\cite{owen2005quasi}, the Monte Carlo schemes do not easily extend to
low discrepancy sequences.\\
Caflisch~\cite{caflisch98} already notes that the simplest way of sampling from both pseudo- and quasi-random numbers is inverse transform sampling (ITS) using the
inverse cumulative probability density. This method, especially for higher dimensions is also known as Rosenblatt-Mück transformation~\cite{rosenblatt1952remarks,Hlawka1972}.\\
For a one dimensional probability density $p_{\mathbf{X}} : [x_{\min},x_{\max}] \rightarrow \mathbb{R}$ the corresponding cumulative density reads
\begin{equation}
 P_{\mathbf{X}} (x) = \int_{x_{\min}}^{x}  p(x) \mathrm{d} x.
\end{equation}
Given a uniformly distributed pseudo- or quasi-random number $\mathbf{u} \sim \mathcal{U}(0,1)$ the corresponding sample $\mathbf{x}$ 
from the probability density $p_{\mathbf{X}}$ is obtained by using the inverse cumulative distribution function $ P^{-1}_{\mathbf{X}} $
\begin{equation}
 P_{\mathbf{X}} (\mathbf{x}) = \mathbf{u}  \quad \Rightarrow \quad \mathbf{x} =  P^{-1}_{\mathbf{X}} (\mathbf{u}).
\end{equation}
Note that it is also common to solve the inversion by Picard iterations or a Newton method.\\
For a given initial condition $f(x,v,t=0)$ to be used in a PIC code the density is mostly so simple that
it decomposes into a tensor product of one dimensional pieces~\cite{caflisch98}, which are then sampled by the one dimensional ITS. \\
We know that PIC codes perform poorly in situations with small perturbations, which is mostly the case in the initial phase of a simulations. Therefore, one might
use a spectral solver to start the simulation for $t \in [0, t_0]$  and then continue with PIC for $t \in [t_0, t_{\max}]$ by using the density
$f(x,v,t=t_0)$ as the initial condition. This requires importance sampling from the density $f(x,v,t=t_0)$ 
which cannot be done anymore by ITS one dimensional pieces. Therefore, we introduce two dimensional (inverse) transform sampling~\cite{olver2013fast} for sampling from
an arbitrary PDF $ g_{\mathbf{X},\mathbf{V}}(x,v)$.\\
On starts with a sample in the first dimension, which is obtained by ITS from a marginal density. Then as we walk through the dimensions we inverse transform sample
from the conditional marginal distribution (integrating over all the higher dimension)
\textit{given} all the previous samples from the lower dimensions.\\
In the two dimensional case we are given a uniformly distributed pseudo- or quasi-random tuple $(\mathbf{u}_x,\mathbf{u}_v) \in \mathcal{U} (0,1)^2$ and have to obtain
the sample $(\mathbf{x},\mathbf{v})$. The first marginal distribution reads
\begin{equation}
 g_{\mathbf{X}} (x) = \int_{ v_{\min}}^{v_{\max}}  g_{\mathbf{X},\mathbf{V}}(x,v)~  \mathrm{d}v,
\end{equation}
and can be sampled from by finding $ \mathbf{x}$ such that $G_{\mathbf{X}} (\mathbf{x})=\mathbf{u}_{x}$, where
\begin{equation}
  G_{\mathbf{X}} (x) =  \int_{ v_{\min}}^{v_{\max}} \int_0^x  g_{\mathbf{X},\mathbf{V}}( \hat{x},v)~  \mathrm{d}\hat{x}\mathrm{d}v
\end{equation}
This is nothing else than using $\mathbf{u}_x$ for ITS from the marginal density $g_{\mathbf{X}}$.
Now \textit{given} the sample $\mathbf{x}$ the conditional density for the likelihood of having a particle at $v$ reads
\begin{equation}
 g_{\mathbf{X}= \mathbf{x}, \mathbf{V}} (v) = \frac{  g_{\mathbf{X},\mathbf{V}}(\mathbf{x},v)  }{ g_{\mathbf{X}} (\mathbf{x})  } 
 = 
 \frac{  g_{\mathbf{X},\mathbf{V}}(\mathbf{x},v)  }{ \int_{ v_{\min}}^{v_{\max}}  g_{\mathbf{X},\mathbf{V}}(\mathbf{x},v)~  \mathrm{d}v  }  
\end{equation}
Note that $g_{\mathbf{X},\mathbf{V}}$ gives as the probability of having both $(x,v)$, but since we already fixed $x=\mathbf{x}$ we have to normalize 
with the corresponding probability for $\mathbf{x}$, namely $g_{\mathbf{X}}(\mathbf{x})$.
After inverting the corresponding cumulative conditional probability density,
\begin{equation}
  G_{\mathbf{X}= \mathbf{x}, \mathbf{V}} (v) =  \int_{ v_{\min}}^{v}  ~ g_{\mathbf{X}= \mathbf{x}, \mathbf{V}} ( \hat{v})~\mathrm{d}\hat{v}
  =
    \frac{  \int_{ v_{\min}}^{v}  g_{\mathbf{X},\mathbf{V}}(\mathbf{x},\hat{v})  ~\mathrm{d}\hat{v} }{ \int_{ v_{\min}}^{v_{\max}}  g_{\mathbf{X},\mathbf{V}}(\mathbf{x},\hat{v})~  \mathrm{d}\hat{v}  }  
\end{equation}
according to
\begin{equation}
G_{\mathbf{X}= \mathbf{x}, \mathbf{V}} (\mathbf{v}) = \mathbf{u}_v
\end{equation}
the second sample $\mathbf{v}$ is obtained. Considering the following map
\begin{equation}
 \begin{split}
\Pi:~& [x_{\min},x_{\max}] \times [v_{\min},v_{\max}] \rightarrow [0,1]^2\\
  &(x,v) \mapsto 
 \begin{pmatrix}
   G_\mathbf{X}( x) \\     G_{\mathbf{X}= x, \mathbf{V}} (v )   
 \end{pmatrix}
 =
\begin{pmatrix}
    \int_{ v_{\min}}^{v_{\max}} \int_0^x  g_{\mathbf{X},\mathbf{V}}( \hat{x},\hat{v})~  \mathrm{d}\hat{x}\mathrm{d} \hat{v} \\    
    \frac{  \int_{ v_{\min}}^{v}  g_{\mathbf{X},\mathbf{V}}(x,\hat{v})  ~\mathrm{d}\hat{v} }{ \int_{ v_{\min}}^{v_{\max}}  g_{\mathbf{X},\mathbf{V}}(x,\hat{v})~  \mathrm{d}\hat{v}  }.
 \end{pmatrix}
 \end{split}
\end{equation}
Then the entire procedure of inverse transform sampling can be described by the inverse map
$\Pi^{-1}:~  [0,1]^2 \rightarrow [x_{\min},x_{\max}] \times [v_{\min},v_{\max}]$. By considering the Jacobi matrix of $\Pi$, 
\begin{multline}
\mathrm{D} \Pi (x,v)=\\
 \begin{pmatrix}
    \int_{ v_{\min}}^{v_{\max}}  g_{\mathbf{X},\mathbf{V}}( x,\hat{v})\mathrm{d} \hat{v}      & 0 \\    
    \frac{ 
\int_{ v_{\min}}^{v}  \partial_x  g_{\mathbf{X},\mathbf{V}}(x ,\hat{v})  \mathrm{d}\hat{v} 
\int_{ v_{\min}}^{v_{\max}} g_{\mathbf{X},\mathbf{V}}(x,\hat{v})  \mathrm{d}\hat{v}    
   -
\int_{ v_{\min}}^{v}  g_{\mathbf{X},\mathbf{V}}(x,\hat{v})  \mathrm{d}\hat{v} 
\int_{ v_{\min}}^{v_{\max}}   \partial_xg_{\mathbf{X},\mathbf{V}}(x,\hat{v})  \mathrm{d}\hat{v}    
}{
    \left( \int_{ v_{\min}}^{v_{\max}}  g_{\mathbf{X},\mathbf{V}}(x,\hat{v})  \mathrm{d}\hat{v} \right)^2 
    } 
    &\frac{  g_{\mathbf{X},\mathbf{V}}(x,v) }{ \int_{ v_{\min}}^{v_{\max}}  g_{\mathbf{X},\mathbf{V}}(x,\hat{v})  \mathrm{d}\hat{v}  } 
 \end{pmatrix}
\end{multline}
the Jacobian introduced by the map $\Pi^{-1}$ is then the Jacobi determinant of $\Pi$:
\begin{equation}
 \mathrm{det} \left(  \mathrm{D} \Pi(x,v)  \right)  =  
 \int_{ v_{\min}}^{v_{\max}}  g_{\mathbf{X},\mathbf{V}}( x,\hat{v})~\mathrm{d} \hat{v}
 \frac{  g_{\mathbf{X},\mathbf{V}}(x,v) }{ \int_{ v_{\min}}^{v_{\max}}  g_{\mathbf{X},\mathbf{V}}(x,\hat{v})~  \mathrm{d}\hat{v}}
 =g_{\mathbf{X},\mathbf{V}}(x,v). 
\end{equation}
This proofs that transforming uniform samples by $\Pi^{-1}$ introduces the Jacobian $g_{\mathbf{X},\mathbf{V}}$ which means we sample from $g_{\mathbf{X},\mathbf{V}}$, which is what we wanted.
The discrete map, based on bilinear interpolation, satisfies smoothness conditions such that the sampling is valid also for QMC numbers, see also~\cite{aistleitner2013low,basu2016transformations}.\\
Suppose a probability density is given in spectral form as
\begin{equation}
\label{fourierg}
  g(x,v) = \sum_{k_x}  \sum_{k_v} \hat{g}(k_x,k_v) \euler^{\imagi ( k_x  (x-x_{\min}) + k_v (v-v_{\min}) )},
\end{equation}
then the marginals can be represented directly by
 \begin{equation}
  G_{\mathbf{X}} (x)  = \hat{g}(0,0) x +   \sum_{k_x\neq 0}  \hat{g}(k_x,0) \frac{1}{ \imagi k_x} \left( \euler^{\imagi  k_x  x-x_{\min})  } -1 \right).
 \end{equation}
Apart from the fact that such dense Fourier interpolation is very expensive another complication arises in the combination with the used spectral solvers.
The obtained distribution function $f$ is not necessarily non-negative such that by the normalization in eqn.~\eqref{normalizef} an absolute value is introduced
in eqn.~\eqref{fourierg} which makes the corresponding anti-derivatives much more complicated. Although it might be physically not reasonable particles can still have the negative weight
$\mathbf{w} = \frac{f( \mathbf{x}, \mathbf{v} ) }{g_{\mathbf{X},\mathbf{V}(\mathbf{x}, \mathbf{v} ) }}$ in order to continue the PIC as close as possible to the spectral solution.
Also, anti-aliasing and Fourier filtering can mitigate the problem but there are no guarantees.
But be warned, ignoring the non-negativity constraint for the sampling density leads to non-monotonic increasing marginal densities
such that the inversion is not well-posed anymore such that this will not work.
Since the Fourier interpolation is expensive anyway and the density $f$ is with the help of the FFT available on a Cartesian grid anyhow, interpolation suggests itself. 
\subsection{Inverse Transform Sampling from a Bilinear Interpolant}
Suppose a two dimensional sampling density $g$ is given as a bilinear interpolant at grid points $(x_i, y_j)$ with values $\left( g_{\mathbf{X,Y}}^{i,j} \right) = g_{\mathbf{X,Y}}(x_i,y_j) $ as
\begin{equation}
\begin{split}
g_{\mathbf{X,Y}}(x,y)= 
\begin{bmatrix}
 1- \frac{ x - x_{i,j}}{\Delta x}  & \frac{ x - x_{i,j}}{\Delta x} 
\end{bmatrix}
\begin{bmatrix}
 g^{i,j} & g^{i,j+1}\\
 g^{i+1,j} & g^{i+1,j+1}\\
\end{bmatrix}
\begin{bmatrix}
 1- \frac{ y - y_{i,j}}{\Delta y}  & \frac{ y - y_{i,j}}{\Delta y} 
\end{bmatrix}\\
\text{ for } x \in [ x_i, x_{i+1}] \text{ and } y \in [ y_j, y_{j+1}]
\end{split}.
\end{equation}
Then the integral over the entire domain is given by the trapezoidal rule
\begin{equation} 
  \int_{y_{\min}}^{y_{\max}} \int_{x_{\min}}^{x_{\max}}  g_{\mathbf{X,Y}}(x,y) ~\mathrm{d}x \mathrm{d}y =
  \sum_{i=1}^{N_x-1}  \sum_{j=1}^{N_y-1} \left(g^{i,j} + g^{i+1,j} +g^{i,j+1} +g^{i+1,j+1} \right) \frac{\Delta x \Delta y}{4} = 1.
\end{equation}
In the case that $(g_{\mathbf{X,Y}}^{i,j})\geq 0$ for all $i,j$ the bilinear interpolation guarantees the positivity of the interpolant $g_{\mathbf{X,Y}}$.
Provided a pair of uniform random or quasi-random numbers $ (u_x,u_y) \in [0,1]^2$ we describe the step by step procedure for determining the unique sample $(\mathbf{x},\mathbf{y})$
that corresponds to inverse transform sampling of the density according by inverse transform sampling of the density $g_{\mathbf{X,Y}}$.
For sampling in the first dimension we need the marginal density
\begin{equation}
 g_{\mathbf{X}} (x)  = \begin{bmatrix}
 1- \frac{ x - x_{i}}{\Delta x}  & \frac{ x - x_{i}}{\Delta x} 
\end{bmatrix}
\begin{bmatrix}
g_{\mathbf{X}}^i \\ g_{\mathbf{X}}^{i+1}
\end{bmatrix}
 \text{ for } x \in [ x_i, x_{i+1}]
\end{equation}
where the $g_{\mathbf{X}}$ at the grid points is exactly given by the trapezoidal rule as
\begin{equation}
g_{\mathbf{X}}^i:= g_{\mathbf{X}}(x_i) =  \int_{y_{\min}}^{y_{\max}}   g_{\mathbf{X,Y}}(x_i,y) ~\mathrm{d}y 
=   \sum_{j=1}^{N_y-1} \frac{  \left( g^{i,j+1} +g^{i,j+1} \right)\Delta y}{2}.
\end{equation}
This allows us to calculate the cumulative distribution function of the marginal density in $\mathbf{X}$ as
\begin{equation}
  G_{\mathbf{X}} (x) =  \Delta x              
  \begin{bmatrix}
\frac{ x - x_{i}}{\Delta x} (  1- \frac{ x - x_{i}}{2 \Delta x} )  &  \left(\frac{ x - x_{i}}{\Delta x} \right)^2
\end{bmatrix}
\begin{bmatrix}
g_{\mathbf{X}}^i \\ g_{\mathbf{X}}^{i+1}
\end{bmatrix}
+ \sum_{1 \leq k<i } g_{\mathbf{X}}^i 
\end{equation}
Given a uniform random number $ \mathbf{u}_x \in [0,1]$, we search for $\mathbf{x}$ such that
\begin{equation}
 G_{\mathbf{X}} ( \mathbf{x} ) = \mathbf{u}_x.  
\end{equation}
Since the underlying polynomial is only of quadratic type and monotonic increasing this inversion can be solved directly. 
The cell index $i$ of $\mathbf{x}$ is found by determining the largest $i$ such that
\begin{equation}
 \sum_{1 \leq k<i } g_{\mathbf{X}}^i  < \mathbf{u}_x,
\end{equation}
which gives
\begin{equation}
\label{eqn:bilinear:cellinversionx}
 \mathbf{x} = x_i + \Delta x
 \frac{ - g_{\mathbf{X}}^i + \sqrt{ ( g_{\mathbf{X}}^i )^2 + 2 (g_{\mathbf{X}}^{i+1} - g_{\mathbf{X}}^i ) 
\left( \mathbf{u}_x - \sum_{1 \leq k<i } g_{\mathbf{X}}^i \right) \frac{1}{\Delta x} }  }{ g_{\mathbf{X}}^{i+1} - g_{\mathbf{X}}^i }.
\end{equation}
Recall that for a given $ \mathbf{x} \in [ x_i, x_{i+1}]$ the conditional distribution function along the second axis is then obtained by 
\begin{equation}
 g_{\mathbf{X}=\mathbf{x},\mathbf{Y}}(y) = \frac{g_{\mathbf{X,Y}}}{\int_{ y_{\min}}^{y_{\max}} g_{\mathbf{X,Y}}(\mathbf{x},y)~\mathrm{d}y} 
 =  \frac{g_{\mathbf{X,Y}}( \mathbf{x} ,y) }{g_{\mathbf{X}}( \mathbf{x}) },
\end{equation}
where $g_{\mathbf{X}}( \mathbf{x})$ is merely a normalization.
The cumulative conditional distribution function along the second axis reads 
\begin{multline}
G_{\mathbf{X}=\mathbf{x},\mathbf{Y}}(y) =   \Bigg\{
\begin{bmatrix}
 1- \frac{ \mathbf{x} - x_{i,j}}{\Delta x}  & \frac{ \mathbf{x} - x_{i,j}}{\Delta x} 
\end{bmatrix}
\begin{bmatrix}
 g^{i,j}  & g^{i,j+1}\\
 g^{i+1,j} & g^{i+1,j+1}\\
\end{bmatrix}
\begin{bmatrix}
 \frac{ y - y_{i,j}}{\Delta y}  - \frac{1}{2} \left( \frac{ y - y_{i,j}}{\Delta y} \right)^2 & \frac{1}{2} \left( \frac{ y - y_{i,j}}{\Delta y} \right)^2  
\end{bmatrix} \Delta y\\
+ 
\begin{bmatrix}
 1- \frac{ \mathbf{x} - x_{i,j}}{\Delta x}  & \frac{ \mathbf{x} - x_{i,j}}{\Delta x} 
\end{bmatrix}
\begin{bmatrix}
  \sum_{1 \leq k < j}  (g^{i,k}+ g^{i,k+1} ) \frac{\Delta y}{2}\\
  \sum_{1 \leq k < j}  (g^{i+1,k}+ g^{i+1,k+1}) \frac{\Delta y}{2}
\end{bmatrix} \Bigg\} \frac{1}{g_{\mathbf{X}}( \mathbf{x})}
\end{multline}
Once again, given the second uniform random number $ \mathbf{u}_y \in [0,1]$ and the first sample $\mathbf{x}$, we search for the sample $\mathbf{y}$ in the second direction such that
\begin{equation}
G_{\mathbf{X}=\mathbf{x},\mathbf{Y}}(\mathbf{y}) = \mathbf{u}_y.  
\end{equation}
Here the cell index $j$ of $\mathbf{y}$ is found by determining the largest $j$ such that
\begin{equation}
\sum_{1 \leq k < j} \left\{
\begin{bmatrix}
 1- \frac{ \mathbf{x} - x_{i,j}}{\Delta x}  & \frac{ \mathbf{x} - x_{i,j}}{\Delta x} 
\end{bmatrix}
\begin{bmatrix}
    (g^{i,k}+ g^{i,k+1} ) \frac{\Delta y}{2}\\
    (g^{i+1,k}+ g^{i+1,k+1}) \frac{\Delta y}{2}
\end{bmatrix} \right\}
  <  g_{\mathbf{X}}( \mathbf{x}) \mathbf{u}_y,
\end{equation}
where $i$ still denotes the cell index in the first dimension, $\mathbf{x} \in [x_i, x_{i+1}]$.
Modifying eqn.~\eqref{eqn:bilinear:cellinversionx} by scaling the density with the normalization and interpolation in the first dimension the $y$ coordinate reads
\begin{equation}
\begin{split}
 \mathbf{y} &= y_j + \Delta y
 \frac{ - \gamma_0 + \sqrt{ ( \gamma_0 )^2 + 2 (\gamma_1 -\gamma_0 ) 
\left(g_{\mathbf{X}}( \mathbf{x}) \mathbf{u}_x - \delta \right) \frac{1}{\Delta y} }  }{ \gamma_1 -\gamma_0 }\\
\text{with}\\
\gamma_l &= \begin{bmatrix}
        1- \frac{ \mathbf{x} - x_{i,j}}{\Delta x}  & \frac{ \mathbf{x} - x_{i,j}}{\Delta x} 
      \end{bmatrix}
    \begin{bmatrix}
 g^{i,j+l} \\
 g^{i+1,j+l}
\end{bmatrix},
    \quad l=0,1      \\
\delta&= \sum_{1 \leq k < j} \left\{
\begin{bmatrix}
 1- \frac{ \mathbf{x} - x_{i,j}}{\Delta x}  & \frac{ \mathbf{x} - x_{i,j}}{\Delta x} 
\end{bmatrix}
\begin{bmatrix}
    (g^{i,k}+ g^{i,k+1} ) \frac{\Delta y}{2}\\
    (g^{i+1,k}+ g^{i+1,k+1}) \frac{\Delta y}{2}
\end{bmatrix} \right\}.
\end{split}
\end{equation}
\subsection{(Bi)Linear Orthogonal Series Density Estimation}
In order to reconstruct a density from given samples we use the counterpart to the bilinear inverse transform sampling, orthogonal series density estimation with
linear splines. Our linear spline basis $(N^1_j)_{j=1,\dots}$ on a one dimensional grid $(x_j)_{j=1,\dots}$ is defined as
\begin{equation}
\label{linearsplinebasis}
N_i^1(x) = N^1\left(  \frac{x-x_i}{\Delta x} \right), \quad
N^1(x) = \begin{cases}
       1-|x| & \text{ for } |x| \leq 1,\\
           0 & \text{ otherwise}.
          \end{cases}
\end{equation}
Note that the mass matrix required for the $L^2$ projection onto our linear spline basis is sparse has the following coefficients
\begin{equation}
 M^x_{m,n}=  \int N_{m}^1(x) N_{n}^1(x) \dx = \Delta x \begin{cases}
                                    \frac{2}{3}       &\text{ for } m=n,\\
                                    \frac{1}{3}      &\text{ for } |m-n|=1,\\
                                      0    &\text{ otherwise }.
                                  \end{cases} 
\end{equation}
For periodic domains the mass matrix is circulant and can be easily applied as an inverse by the use of the fast Fourier transform~\cite{gray2006toeplitz}.
Here, of course, we consider a two dimensional basis such that the bilinear interpolation can be written with the help of tensor product splines as
\begin{equation}
 g(x,y) = \sum_{i,j} g^{i,j} N_{i,j}(x,y), \quad N_{i,j}(x,y):= N^1\left(  \frac{x-x_i}{\Delta x} \right)N^1\left(  \frac{y-y_j}{\Delta y} \right).
\end{equation}
The corresponding mass matrix for the two dimensional space is obtained as a tensor product $M^{x,y} = M^x \otimes M^y$.
This also means that the $L^2$ projection of a density sampled by bilinear inverse transform sampling onto this bilinear spline space is exact:
\begin{equation}
 \kvec{ g_{i,j} } = \left( M^{x,y} \right)^{-1} \kvec{ \iint g(x,y) N_{i,j}(x,y)  \mathrm{d}x \mathrm{d}{y}          }
\end{equation}
The coefficients $g_{i,j}$ can be approximated as
\begin{equation}
 \kvec{ g_{i,j} } \approx \left( M^{x,y} \right)^{-1} \kvec{ \ \frac{1}{N_s} \sum_{n=1}^{N_s} w_n ~ N_{i,j}(\mathbf{x}_n,\mathbf{y}_n) }
\end{equation}
which corresponds to the cloud in cell scheme combined with the mass matrix. Since we are interested in the sampling density the weights are constant $w_n=1$, but for any other
function they read $w_n= \frac{f(\mathbf{x}_n,\mathbf{y}_n)}{g(\mathbf{x}_n,\mathbf{y}_n)}$.\\
Since we combine this OSDE with spectral densities there is an additional approximation error that we have to review.
Recall that the $m^{\mathrm{th}}$ order B-Spline $N^1(x)$ on a grid of size $h$ is obtained by convolution as
\begin{equation}
\begin{split}
N^m(x) &= \underbrace{ N^0 * \dots * N^0}_{m+1 \text{ times }}(x) = (N^0)^{*(m+1)}(x)=\int_{-\infty}^{\infty} N^0(x-y) N^{m-1}(y) ~\mathrm{d}y\\
\text{ with }
N^0 (x) &= \begin{cases}
        \frac{1}{h} &\text{if } x \in \left(-\frac{h}{2}, \frac{h}{2} \right)\\
        0 & \text{else}.
       \end{cases}\\
\end{split}
\end{equation}
The Fourier transform of one basis function $N^m$ reads
\begin{equation}
 \int_{x_{\min}}^{x_{\max}} N^m(x) \euler^{-i k x} ~\mathrm{d}x  = \left[ \mathrm{sinc} \left( \frac{kh}{2} \right) \right]^{m+1},
\end{equation}
such that we can conclude that the relative error on the $k^{\mathrm{th}}$ Fourier mode is 
\begin{equation}
 \left| 1- \left[ \mathrm{sinc} \left( \frac{kh}{2} \right) \right]^{m+1} \right|
\end{equation}
For our linear splines $m=1$, this can be quite large on the highest mode still represented on the grid $ 1-\mathrm{sinc}(1)^2 \approx  0.2919$, such that we have to decrease the grid size $h$.
Given the Fourier coefficients for the spectral solver this is easily achieved with increasing the resolution by a factor $N_{\mathrm{pad}}$ by zero padding in the inverse FFT.
When applying padding we chose $N_{\mathrm{pad}}=32$, such that the relative error on highest mode reduces to $ 1-\mathrm{sinc}\left( \frac{1}{N_{\mathrm{pad}}} \right)^2 \approx  3.2548e-04$.
\FloatBarrier
\section{Numerical results}
We consider the two dimensional Vlasov--Poisson system with the initial condition
\begin{equation}
\begin{split}
 f(x,v,t=0) &= \frac{1- \epsilon \cos(k  x)}{\sqrt{2\pi}} \left[  (1-n_b)  \euler^{-\frac{v^2}{2}} + \frac{ n_b}{\sigma_b}  \euler^{-\frac{(v-v_b)^2}{2 \sigma_b^2}}      \right] \\
 & \text{ for } x \in [0,L], ~v \in [v_{\min},v_{\max} ], ~ L=\frac{2\pi}{k}\\
\end{split} 
\end{equation}
and the following parameters for two nonlinear test cases:
\begin{center}
\begin{tabular}{c c}
 Landau damping & $ \epsilon= 0.5 , ~ k=0.5, ~[\sigma_b =1 , ~n_b=0]$\\
 Bump-on-tail instability & $ \epsilon= 10^{-3} , ~ k=0.3, ~\sigma_b =0.3 , ~n_b=0.1$\\
\end{tabular}
\end{center}
The number of particles is in the PIC simulation is denoted by $N_p$ and the number of cells by $N_f$. For the Poisson solver finite elements based on cubic B-splines are used.
Random numbers and the quasi random Sobol sequence are provided by MATLAB~\cite{bratley1988algorithm}.
The spectral solver uses the same number of grid points in spatial and velocity space,~$N_x=N_v$.

\FloatBarrier
\subsection{Phase space conservation}
\begin{figure}
\centering
\begin{subfigure}{0.8\textwidth}
 \includegraphics[width=\textwidth]{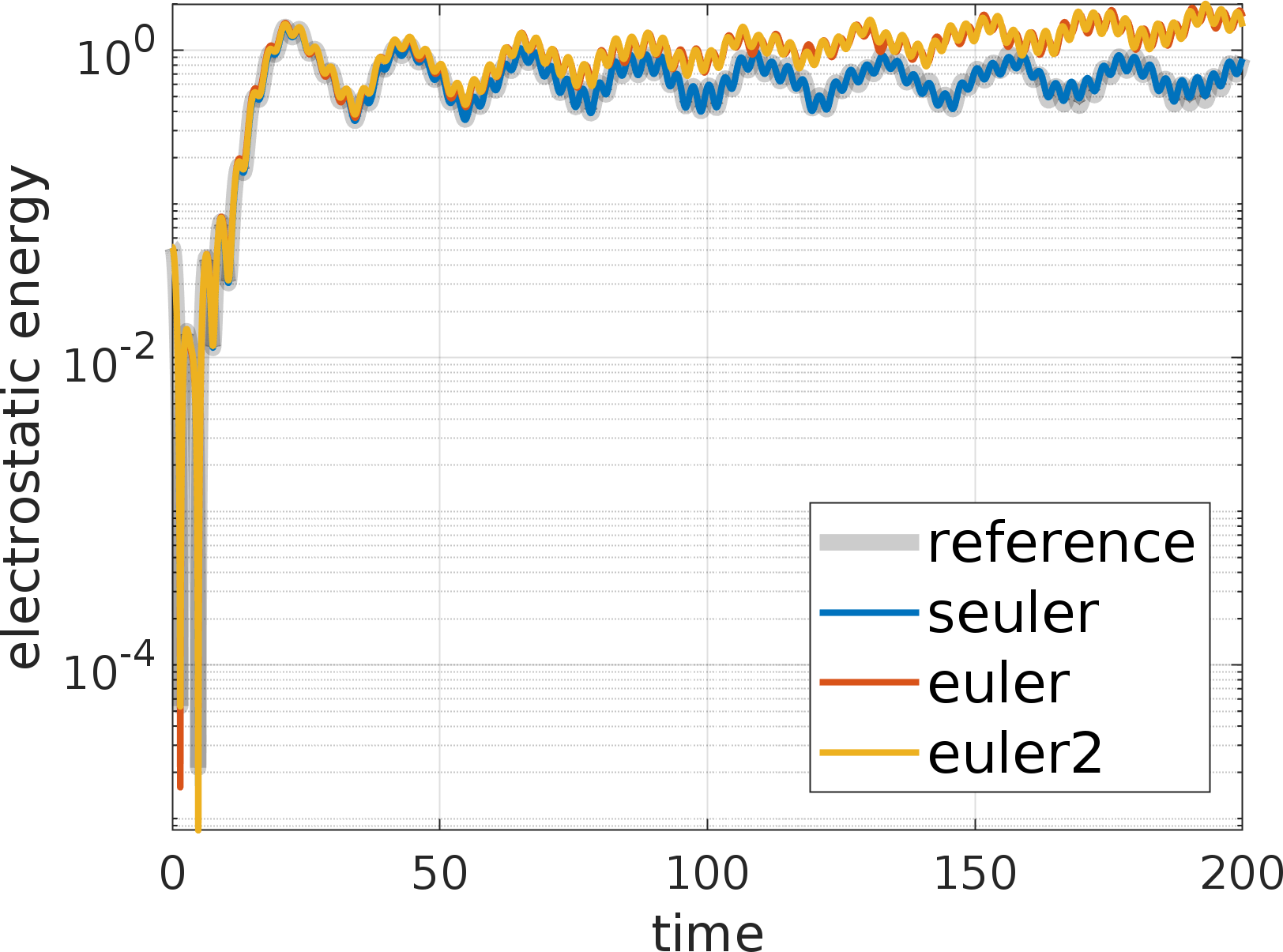} 
 \caption{electrostatic energy}
\end{subfigure}\\
\begin{subfigure}{0.32\textwidth}
\includegraphics[width=\textwidth]{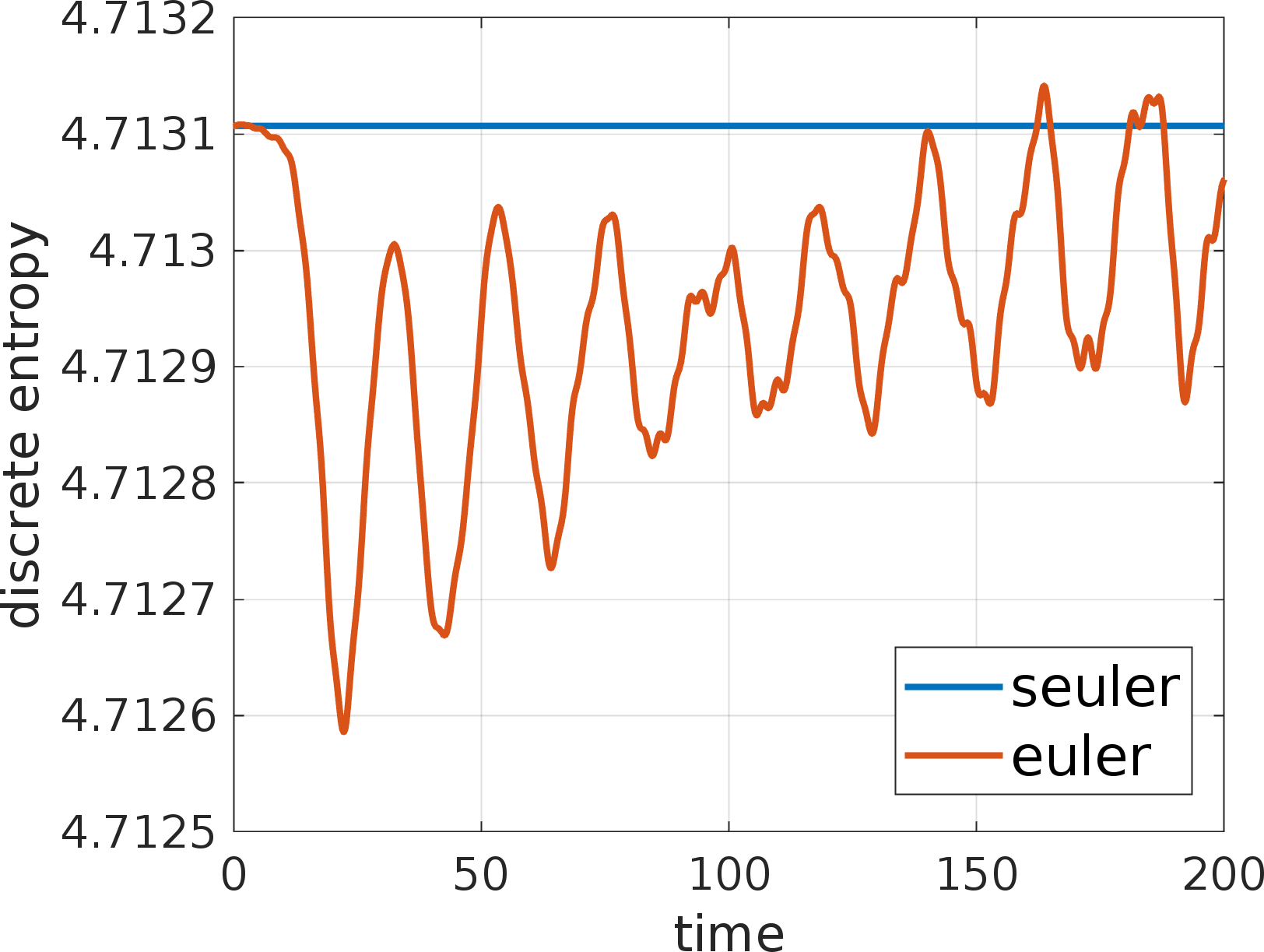}
\caption{discrete entropy}
\end{subfigure}
\begin{subfigure}{0.32\textwidth}
 \includegraphics[width=\textwidth]{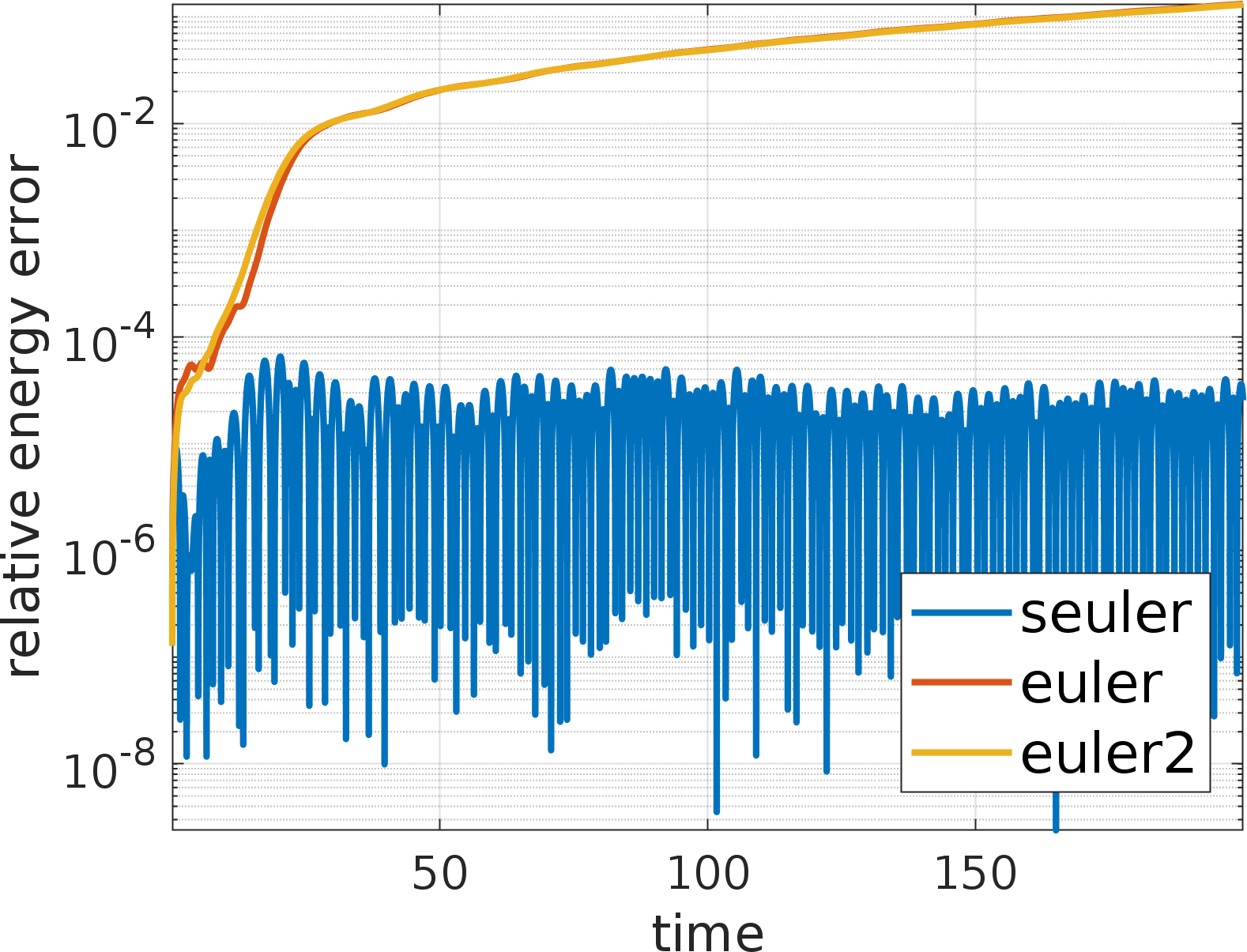} 
 \caption{relative energy error}
\end{subfigure}
\begin{subfigure}{0.32\textwidth}
 \includegraphics[width=\textwidth]{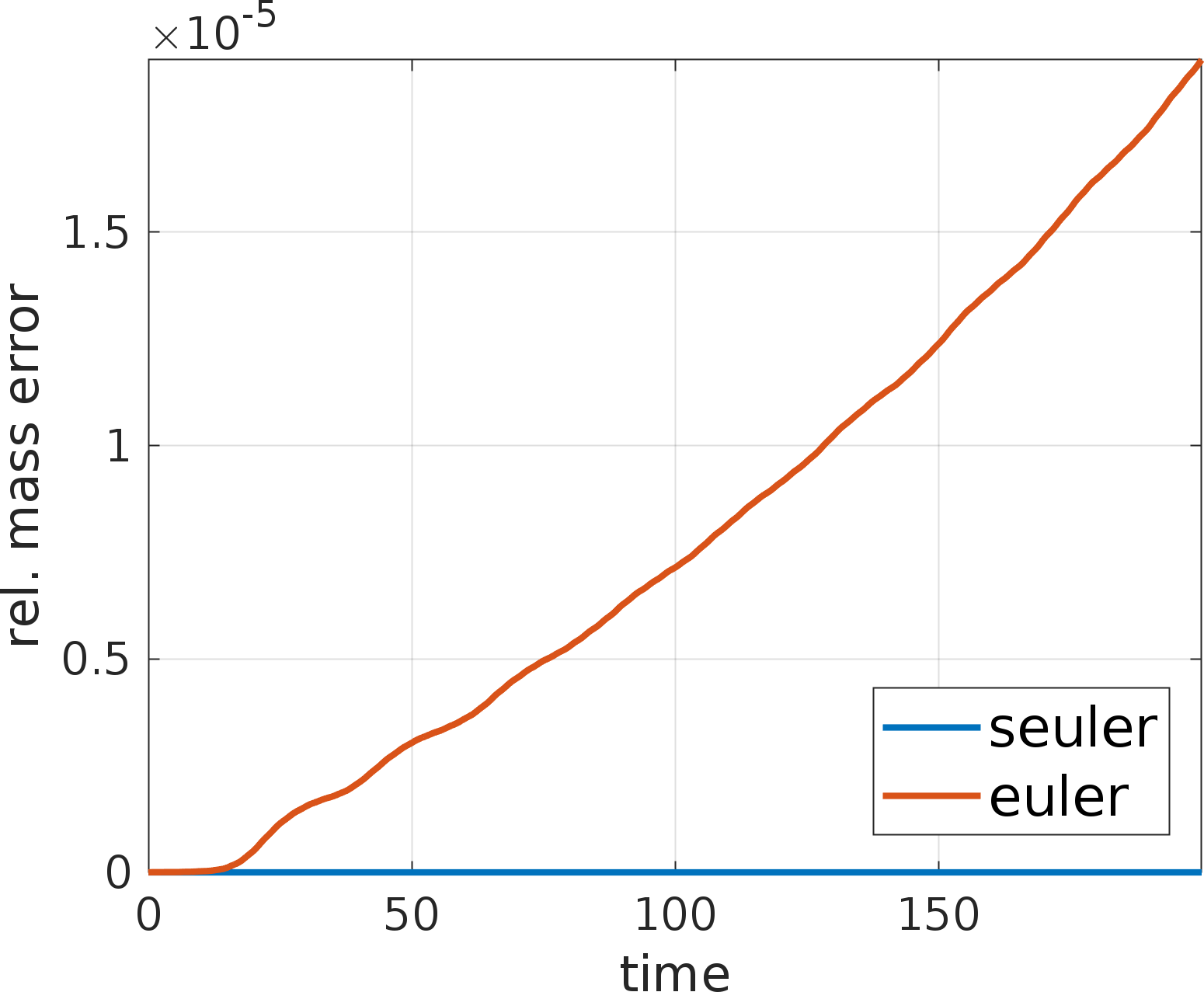} 
 \caption{relative mass error}
\end{subfigure}\\
\caption{PIC simulation using cubic B-Splines $(N_x=32)$ and $\Delta t=0.01$ of a Bump-on-tail instability with $N_p=10^6$
         quasi-randomly distributed particles drawn by inverse transform importance sampling. 
         The explicit Euler (\textit{euler}) is outperformed by the sympletic (\textit{seuler)}).
          }
\label{fig:phasespace_conservation:bumpontail}
\end{figure}

\begin{figure}{H}
\centering
\begin{subfigure}{0.8\textwidth}
 \includegraphics[width=\textwidth]{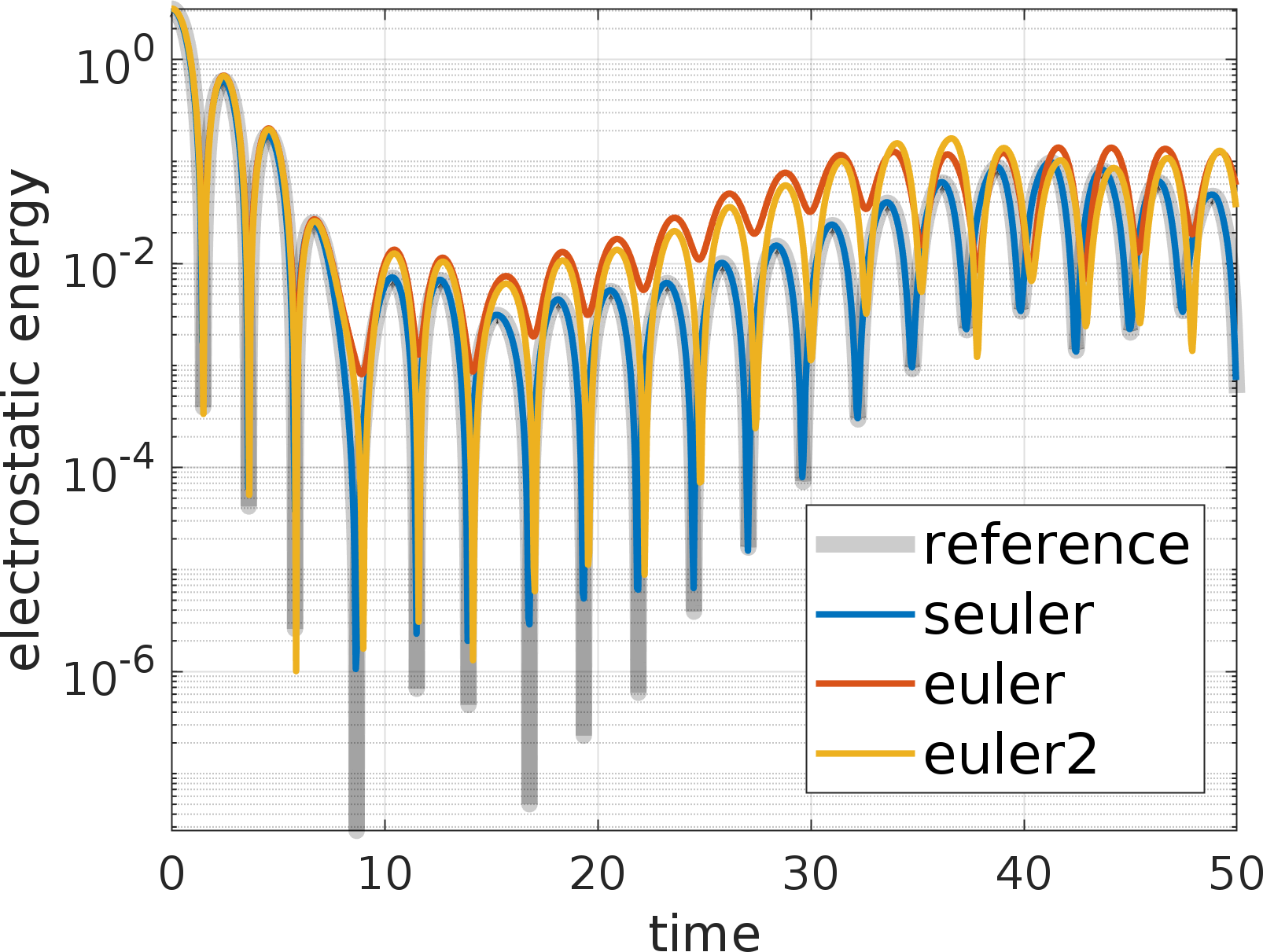} 
 \caption{electrostatic energy}
\end{subfigure}\\
\begin{subfigure}{0.32\textwidth}
\includegraphics[width=\textwidth]{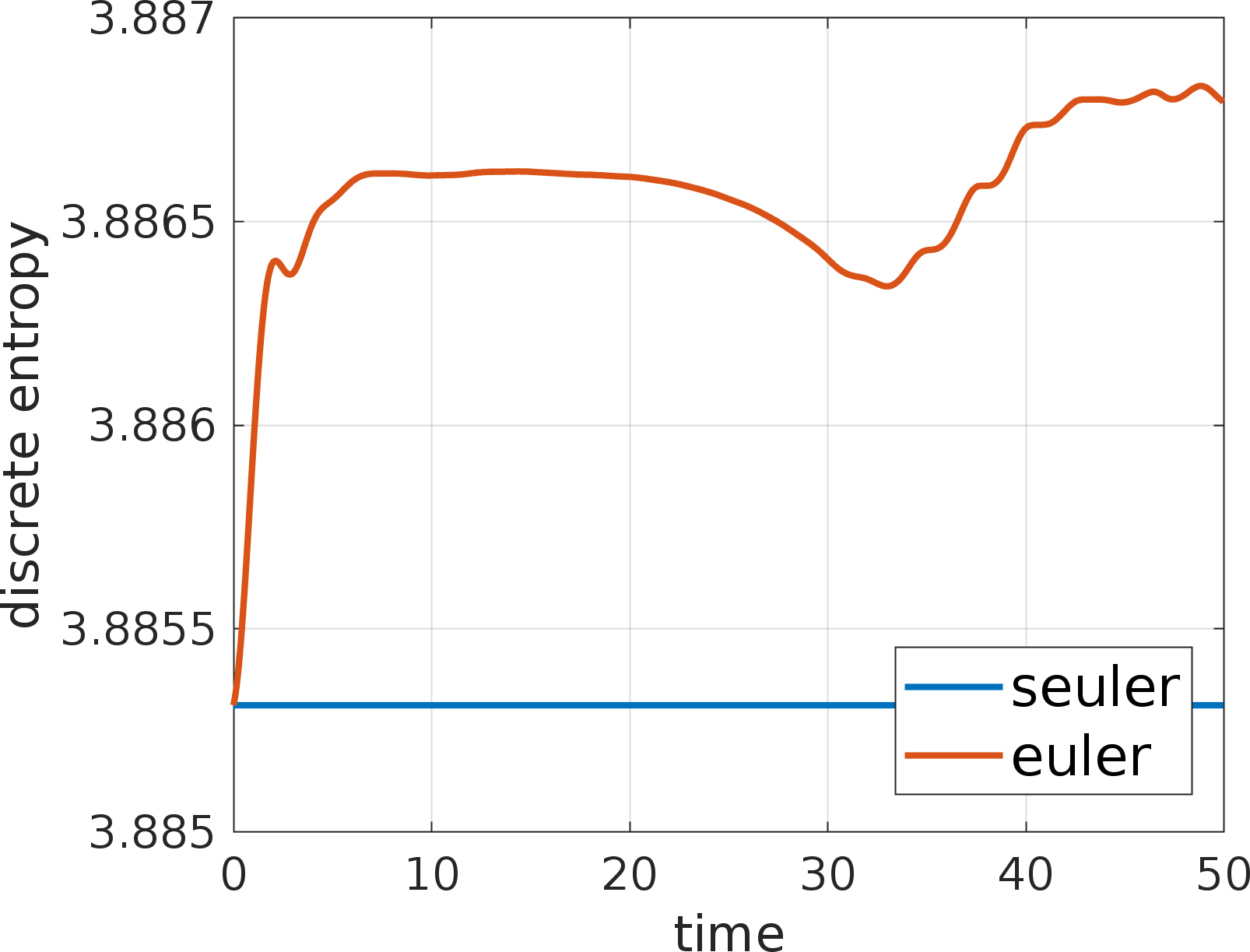}
\caption{discrete entropy}
\end{subfigure}
\begin{subfigure}{0.32\textwidth}
 \includegraphics[width=\textwidth]{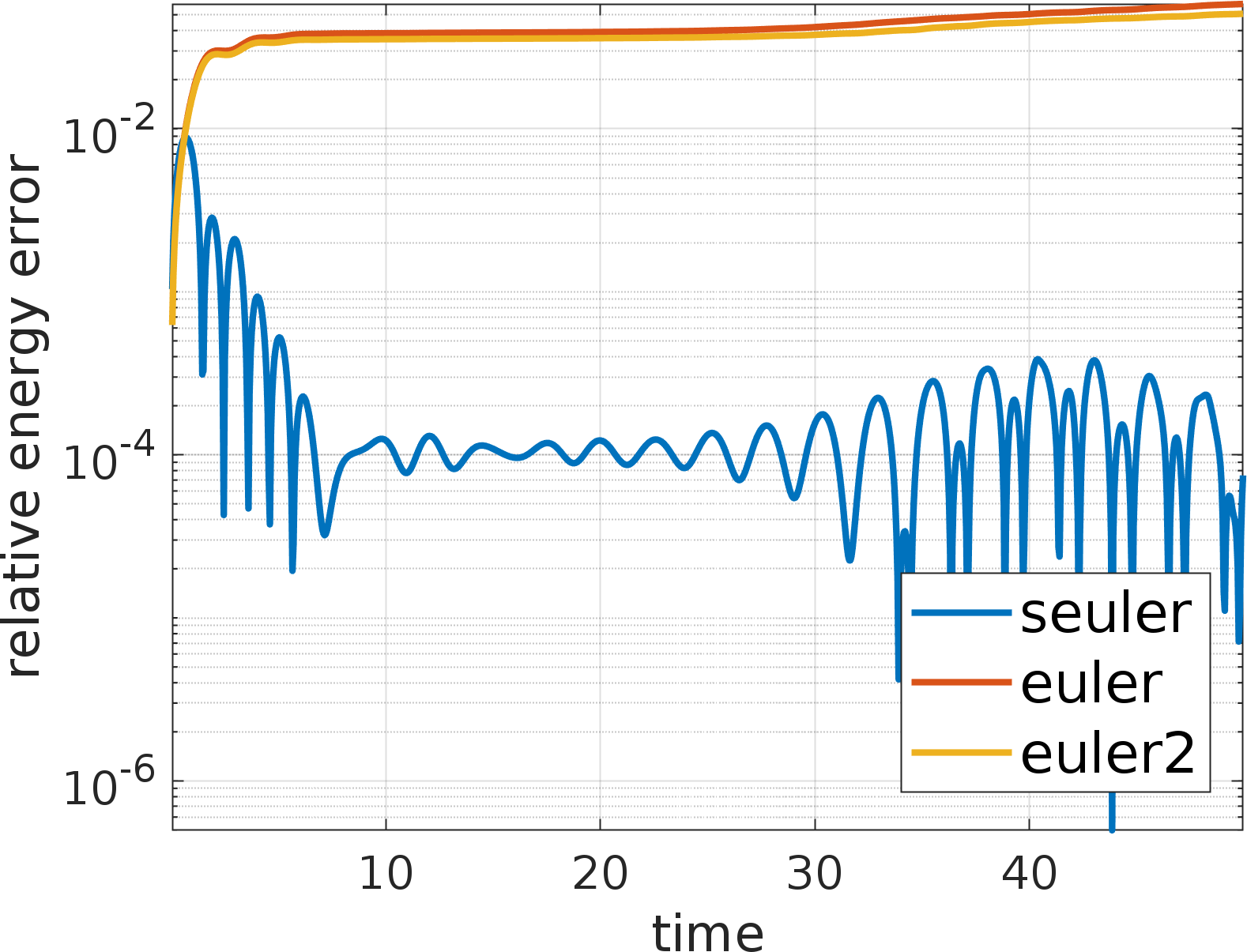} 
 \caption{relative energy error}
\end{subfigure}
\begin{subfigure}{0.32\textwidth}
 \includegraphics[width=\textwidth]{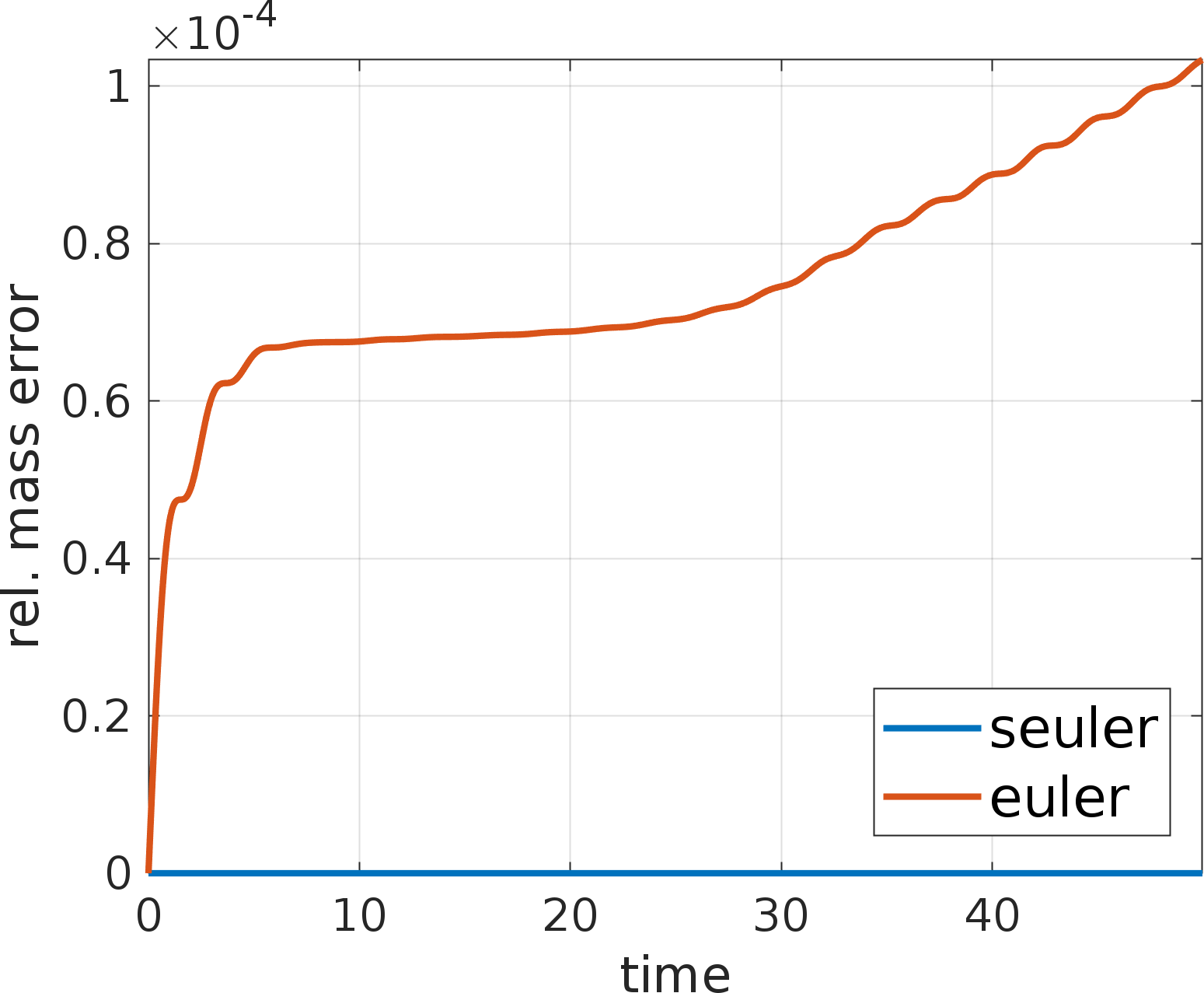} 
 \caption{relative mass error}
\end{subfigure}\\
\caption{PIC simulation of nonlinear strong Landau damping with $N_p=10^6$ quasi-randomly distributed particles drawn by inverse transform importance sampling. 
         Sympletic (\textit{seuler)})
          Cubic B-Spline $(N_x=32)$ }
\label{fig:phasespace_conservation:landau_strong}
\end{figure}

We have already seen in fig.\ref{fig:pic:star_discrepancy} that phase space conserving symplectic Euler as opposed to the explicit Euler preserves the uniformity of a quasi-random sequence implying better convergence.
In general for integrators, which do not preserve phase space volume but are dissipative such as asymptotically preserving schemes 
like~\cite{filbet2016asymptotically,filbet2017asymptotically} the likelihoods have to be propagated accordingly. For the explicit Euler, there are two options.
As we are following the characteristics we leave $f_k$ constant, but rescale the sampling likelihood $g_k$ with the according Jacobi determinant of the flow and call this \textit{euler}.
For \textit{euler2} we ignore the characteristics and also rescale $f_k$ with the Jacobian. Since the likelihoods actually change now it makes sense
to look at otherwise conserved quantities such as the total mass 
\begin{equation}
 \iint f(x,v,t)  \dxdv = \mathbb{E}\left[ \frac{f(\mathbf{X}(t),\mathbf{V}(t),t)}{ g(\mathbf{X}(t),\mathbf{V}(t),t)}  \right] = 
 \frac{1}{N_p} \sum_{k=1}^{N_p} \frac{ f^t_k}{ g^t_k}
\end{equation}
and the discrete variant of the differential entropy 
\begin{multline}
\label{discreteentropy}
  \iint f(x,v,t) \ln(f(x,v,t))  \dxdv\\
  = \mathbb{E}\left[ \frac{f(\mathbf{X}(t),\mathbf{V}(t),t) \ln \left( f(\mathbf{X}(t),\mathbf{V}(t),t)  \right)}{ g(\mathbf{X}(t),\mathbf{V}(t),t)}  \right]
  \approx \frac{1}{N_p} \sum_{k=1}^{N_p} \frac{ f^t_k \ln(f^t_k)}{ g^t_k}.
  \end{multline}
Note, that there are different ways to estimate the entropy from a sample,~\cite{beirlant1997nonparametric} gives an overview, while in ~\cite{ebrahimi1994two}
mesh based examples ready for implementation can be found. By propagating the sampling weight correctly it is possible to observe changes in the entropy caused by the dissipative integrator.
The results can be seen for strong Landau damping in fig.~\ref{fig:phasespace_conservation:landau_strong} and the Bump-on-tail instability in fig.~\ref{fig:phasespace_conservation:bumpontail}.
For the electrostatic field energy a reference solution was calculated using a pseudo spectral solver.

The choices made for \textit{euler} and \textit{euler2} are both inconsistent in some way because the method itself is just not suited for this purpose.
For an asymptotically preserving scheme, where we do not follow the characteristics, it would make sense though to rescale $f_k$ and $g_k$ accordingly.

\FloatBarrier
\subsection{Convergence of OSDE}
\begin{figure}{h}
\centering
\begin{subfigure}{0.45\textwidth}
 \includegraphics[width=\textwidth]{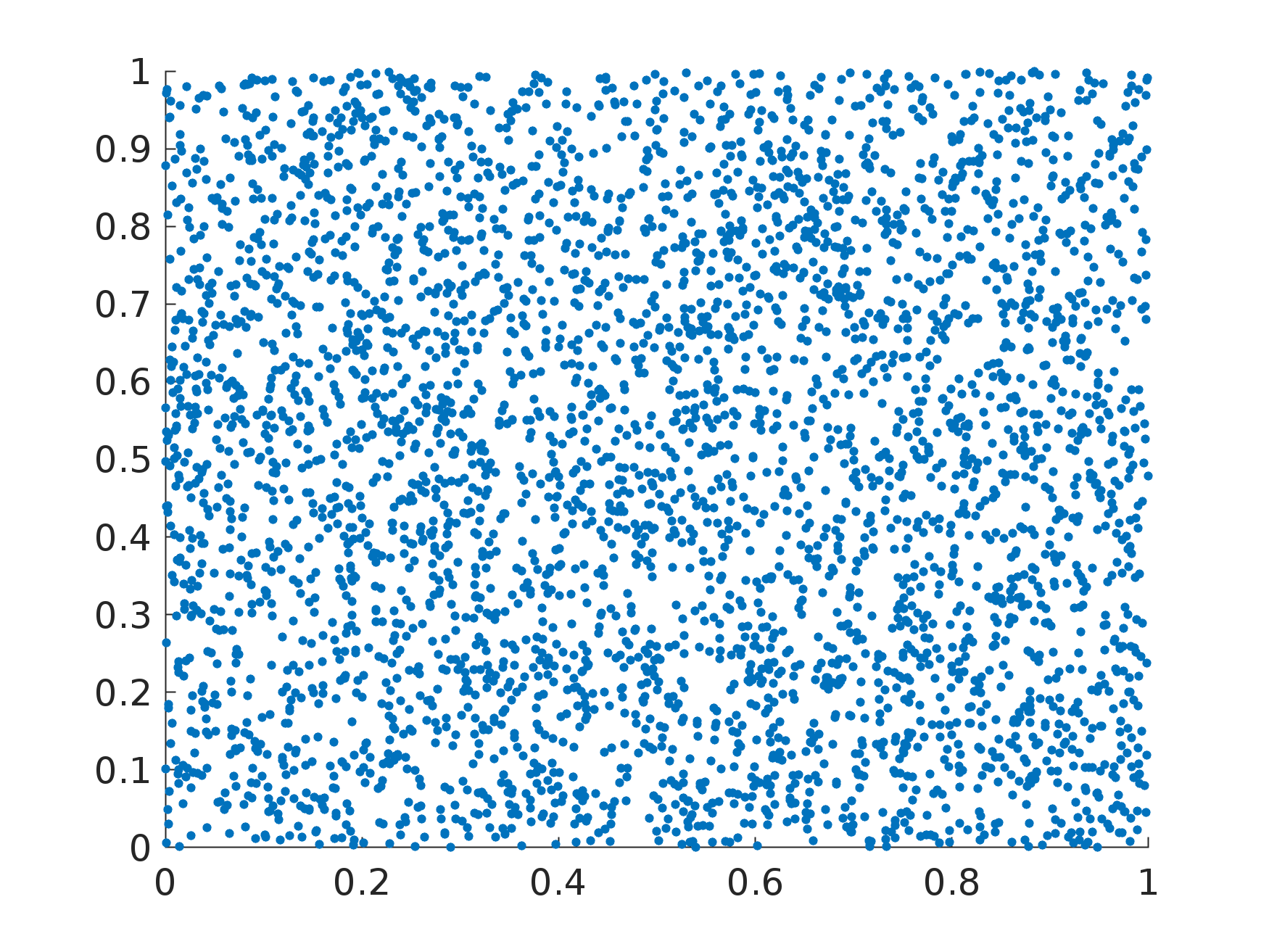} 
 \caption{uniform (pseudo)random numbers}
\end{subfigure}
\begin{subfigure}{0.45\textwidth}
\includegraphics[width=\textwidth]{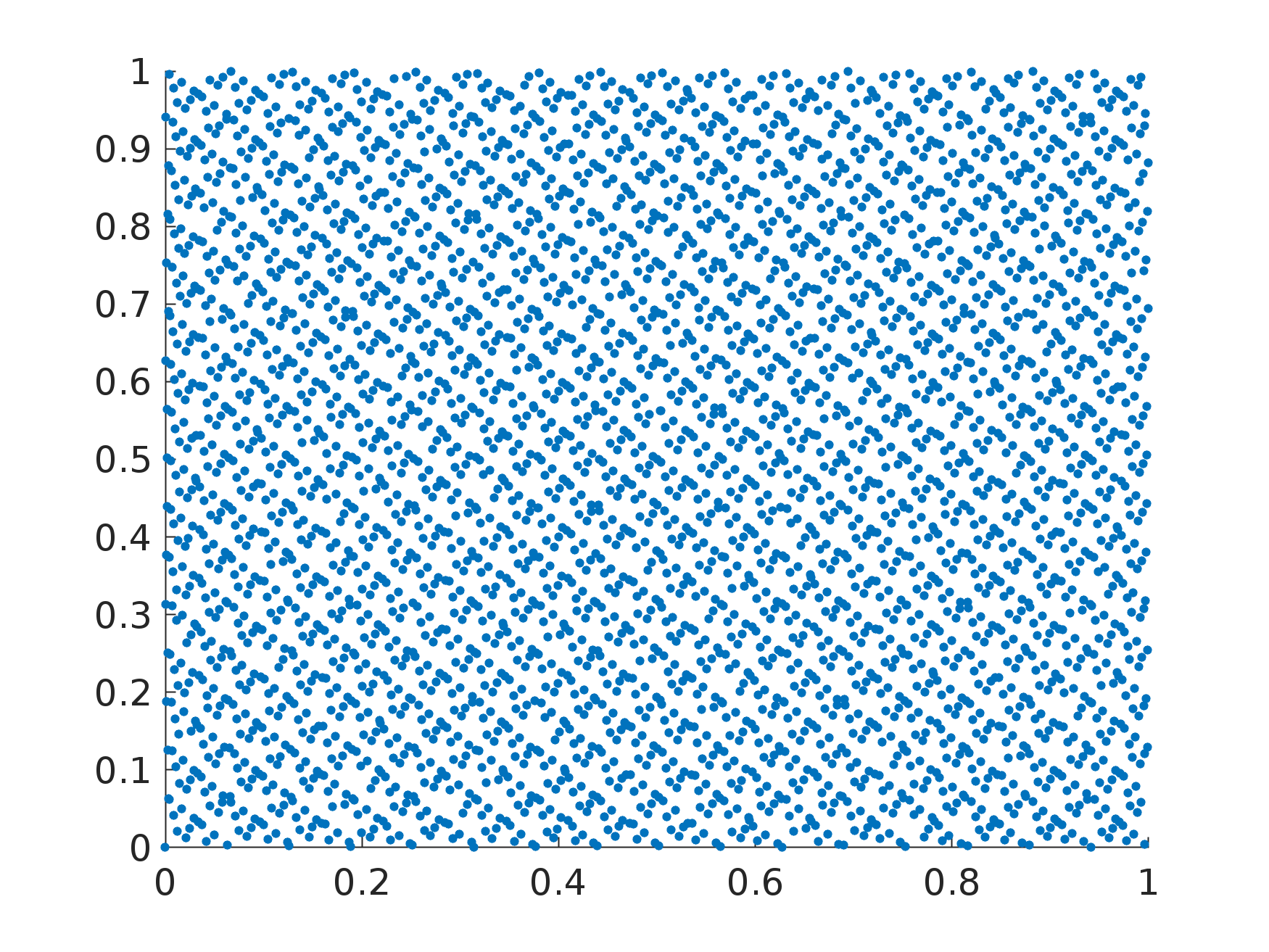}
\caption{uniform Sobol numbers}
\end{subfigure}\\
\begin{subfigure}{0.45\textwidth}
 \includegraphics[width=\textwidth]{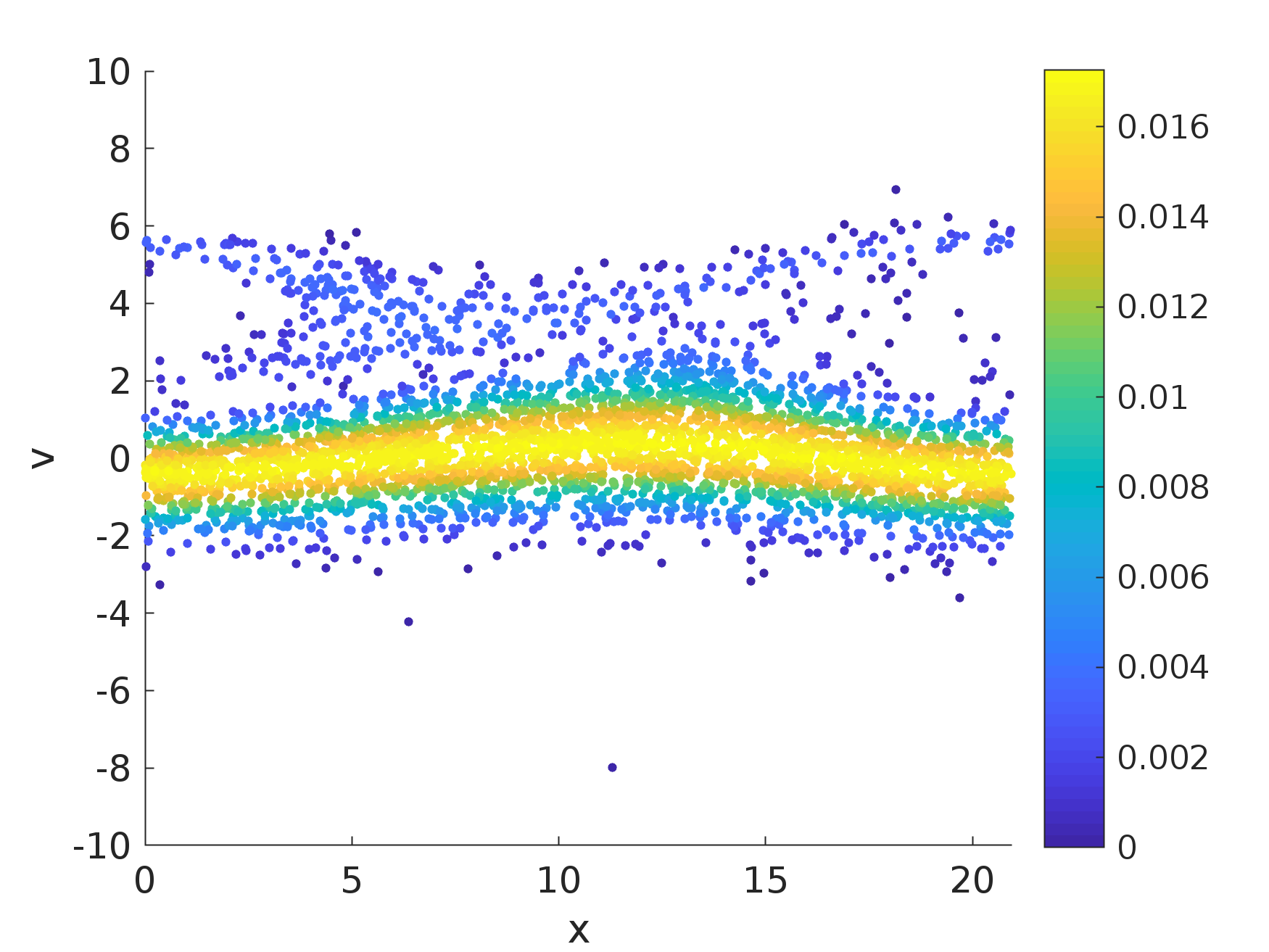} 
 \caption{Random samples}
\end{subfigure}
\begin{subfigure}{0.45\textwidth}
\includegraphics[width=\textwidth]{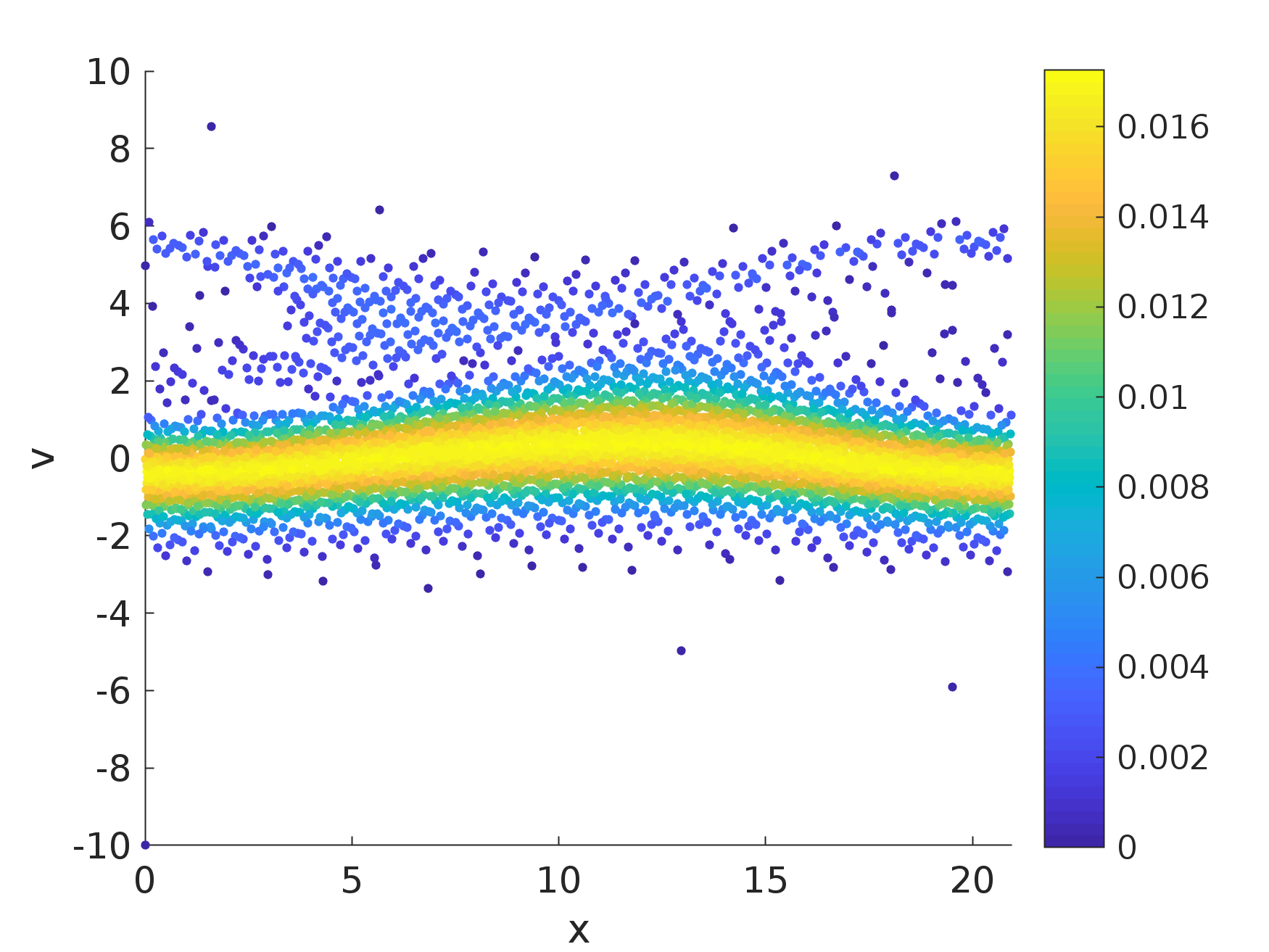}
\caption{Quasi-Random samples}
\end{subfigure}\\
\begin{subfigure}{0.45\textwidth}
\includegraphics[width=\textwidth]{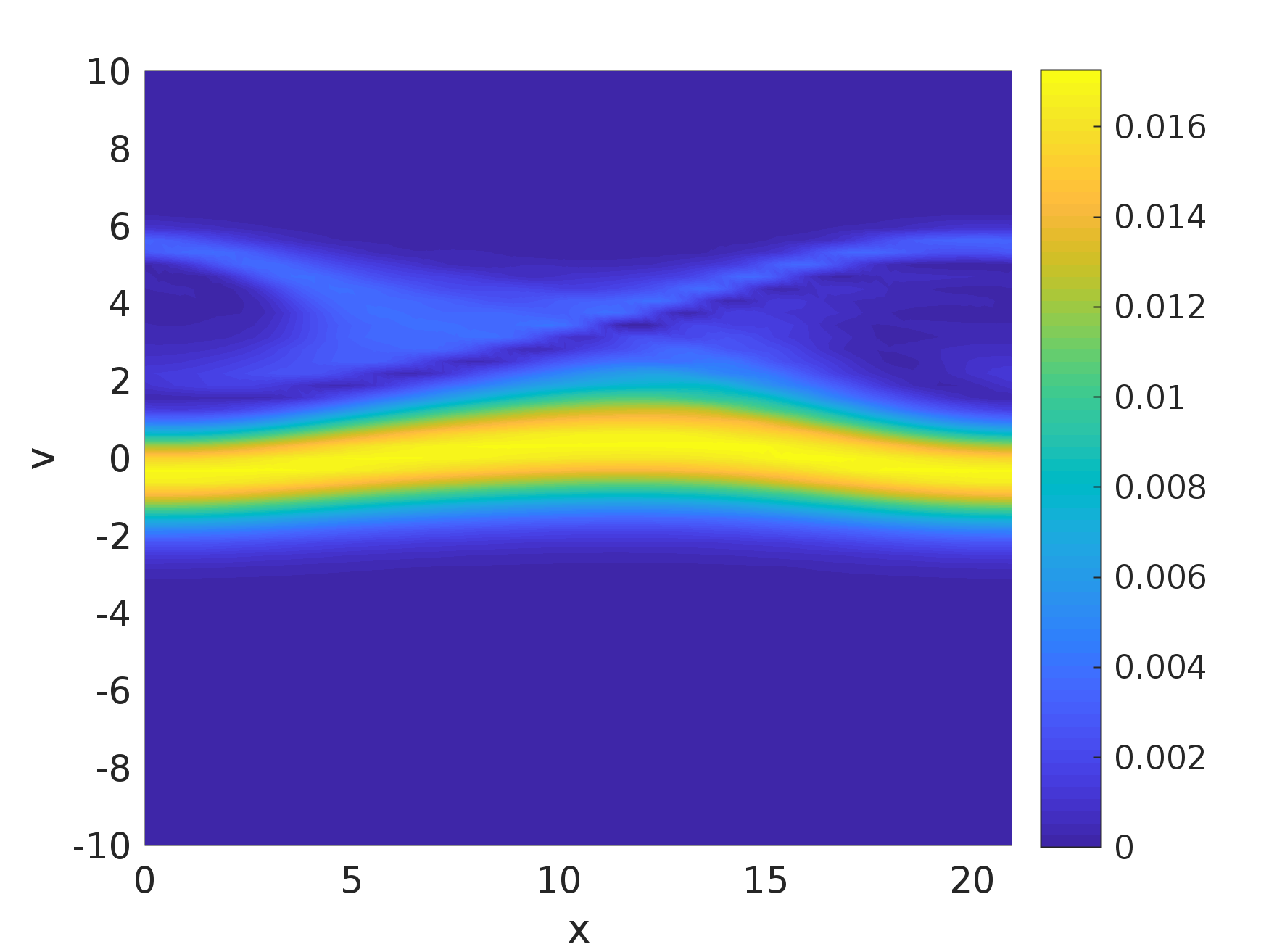}
\caption{sampling density $g$ at $t=30$}
\end{subfigure}
\begin{subfigure}{0.45\textwidth}
\includegraphics[width=\textwidth]{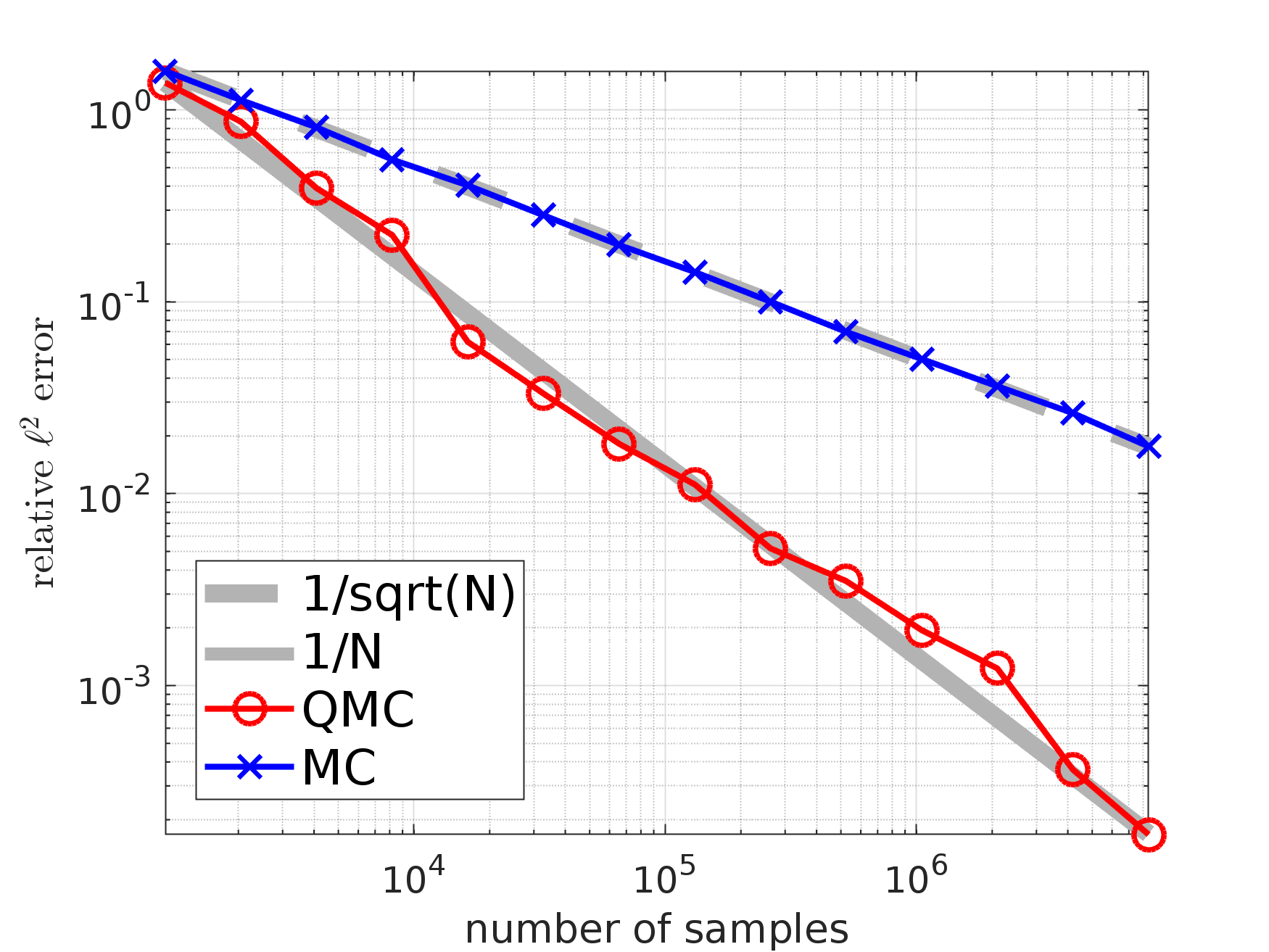}
\caption{Relative $L^2$ sampling error}
\end{subfigure}
\caption{Sampling $N_p=4000$ random and quasi-random particles with bilinear inverse transform sampling from the sampling density $g$ obtained at $t_{\max}=35$ from spectral simulation
          of the Bump-on-tail instability $N_x=N_v=64,~ \Delta t=0.01$.
          Fig. (f) shows the relative $L^2$ error on the density $f$ obtained by OSDE from the samples at increasing number of particles.}
\label{fig:bumpontail_osde}
\end{figure}

\begin{figure}
\centering
\begin{subfigure}{0.45\textwidth}
 \includegraphics[width=\textwidth]{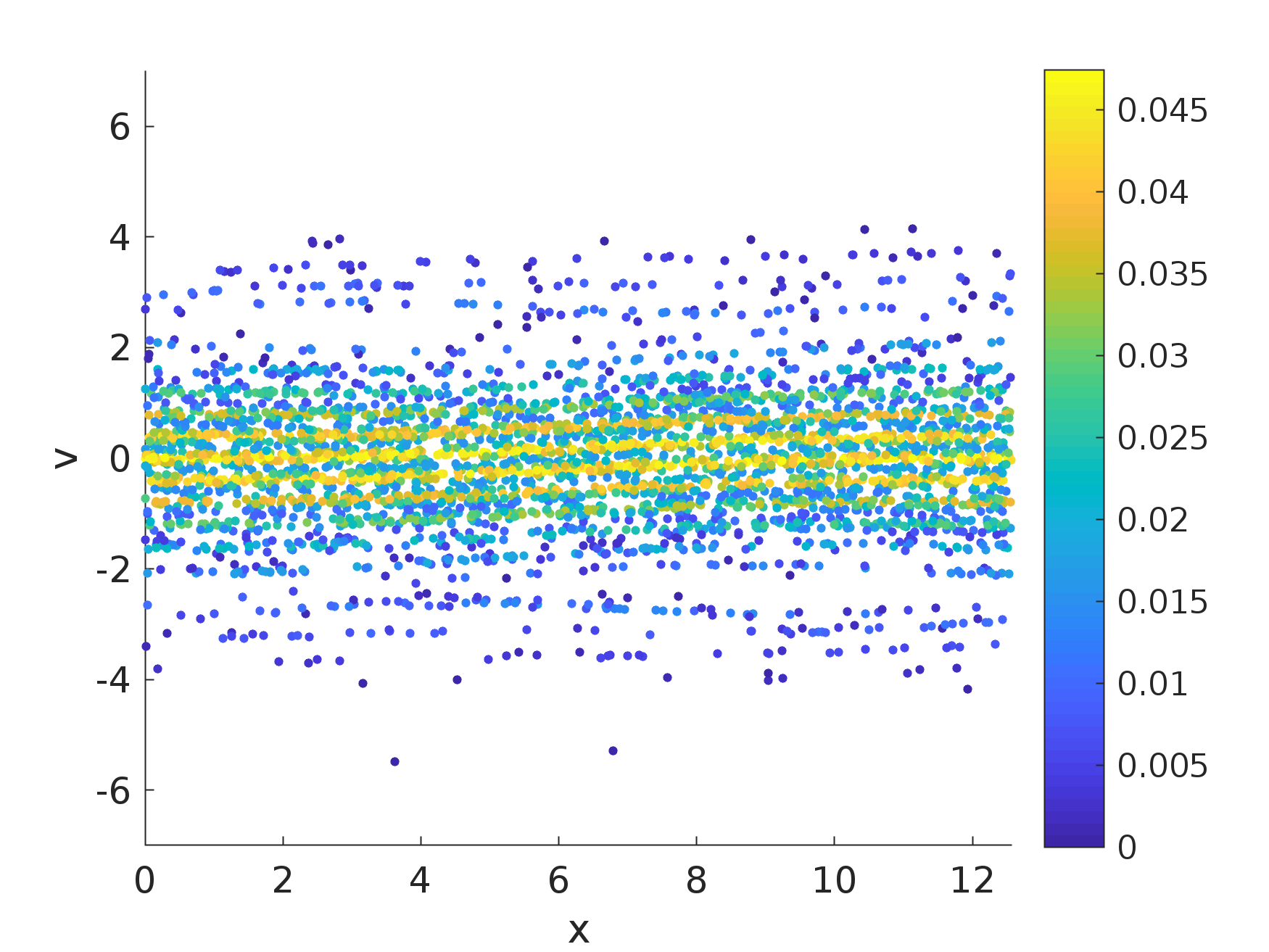} 
  \caption{Random samples}
\end{subfigure}
\begin{subfigure}{0.45\textwidth}
\includegraphics[width=\textwidth]{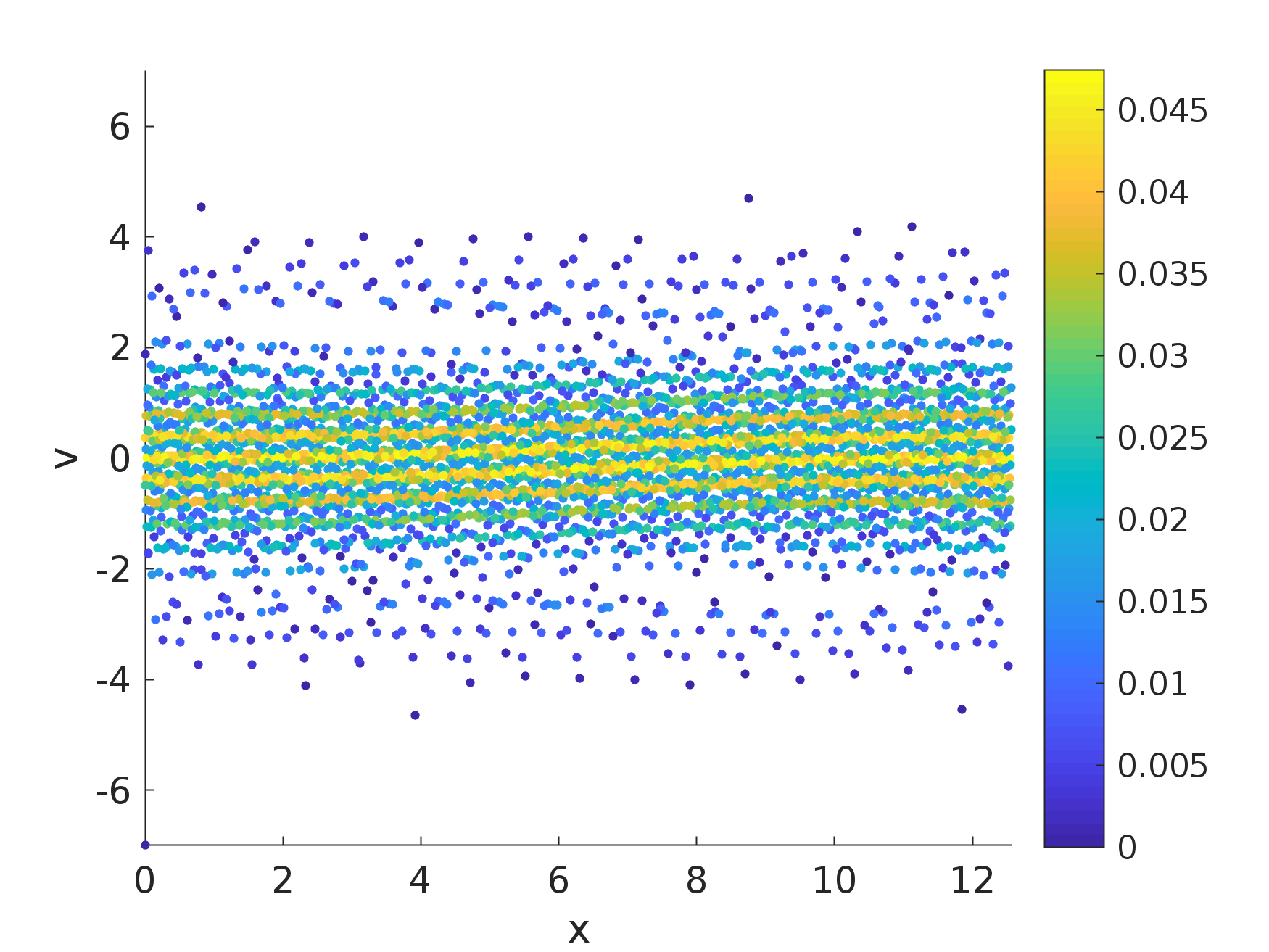}
\caption{Quasi-Random samples}
\end{subfigure}\\
\begin{subfigure}{0.45\textwidth}
\includegraphics[width=\textwidth]{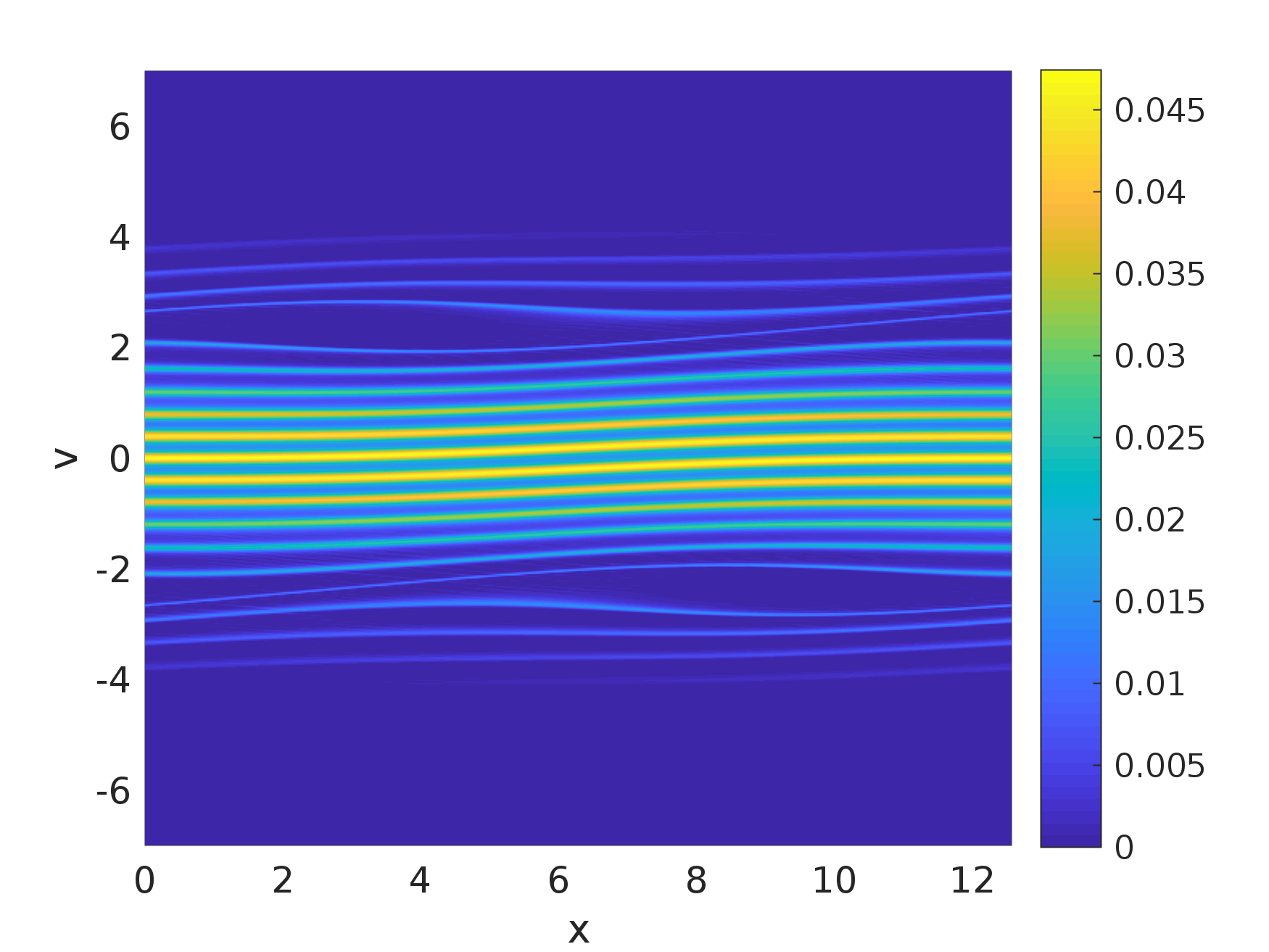}
\caption{sampling density $g$ at $t=30$}
\end{subfigure}
\begin{subfigure}{0.45\textwidth}
\includegraphics[width=\textwidth]{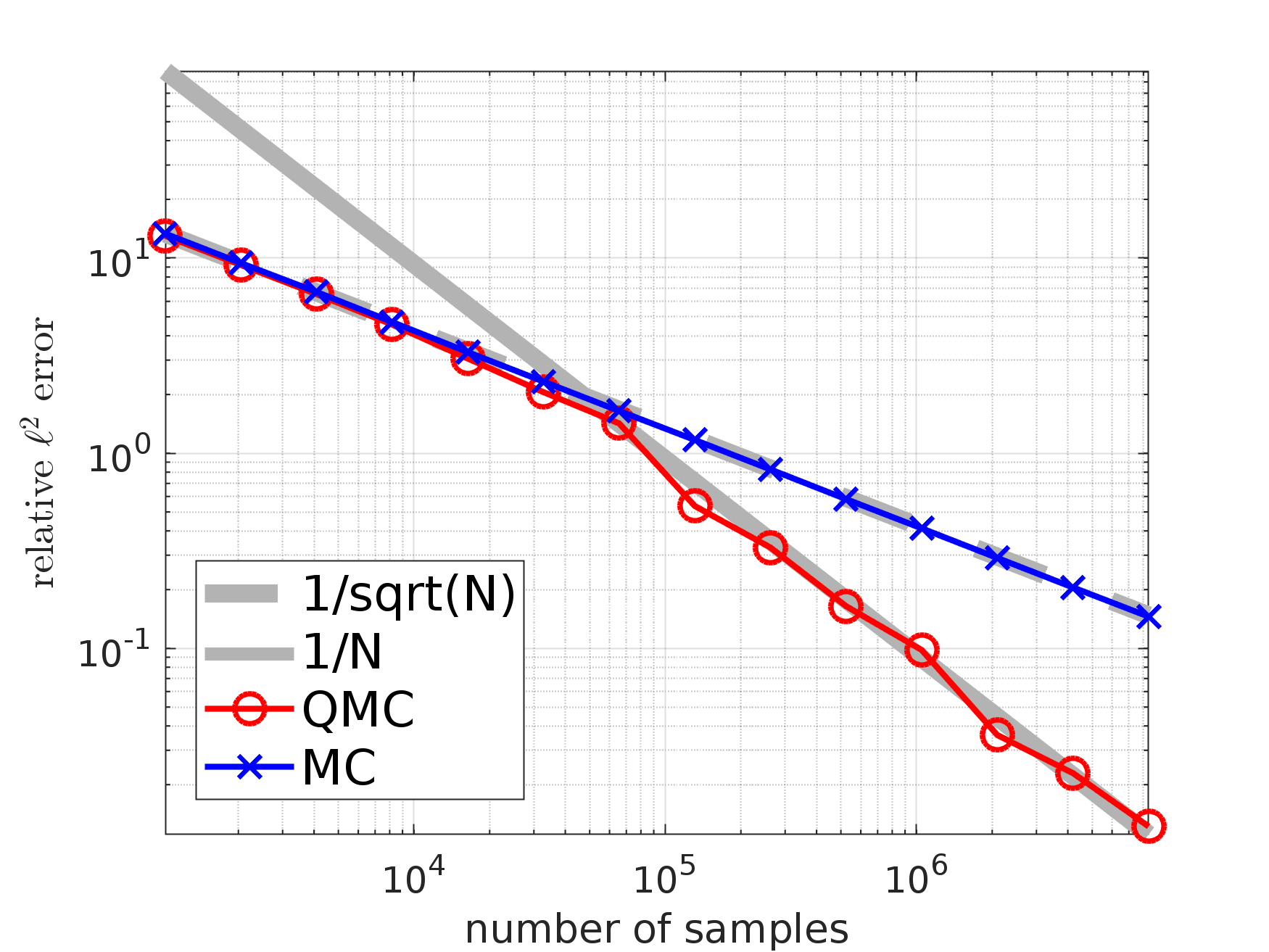}
\caption{Relative $L^2$ sampling error}
\end{subfigure}
\caption{The fine structures of nonlinear Landau damping are resolved at $t_{\max}=30$ by a spectral solver with $N_x=N_v=512,~ \Delta t=0.01$ in (c). The difference between
         $N_p=4000$ random (a) and quasi-random (b) samples obtained from (c) are clearly visible and also seen in (d) in the different convergence rates
         of the corresponding reconstructing OSDE of $f$.}
\label{fig:landau_osde}
\end{figure}

The first experiment should test whether the sampler has been correctly implemented and whether we obtain the expected convergence rates. For this
the Vlasov--Poisson system~\eqref{vlasov1}-\eqref{poisson1} is simulated with the spectral solver for the Bump-on-tail instability and Landau damping in $t \in [0,t_{\max}]$.
The density $f(x,v,t_{\max})$ is obtained on the $N_x \times N_v$ phase-space grid without additional padding. In the following bilinear interpolation
is used for further representation of $f$. Negative values remain and are not truncated.
After obtaining the corresponding sampling density according to eqn.~\eqref{normalizef} samples are drawn by bilinear inverse transform sampling. To verify
that the particles are sampled correctly the density is estimated again by bilinear OSDE.
The convergence rates in fig.~\ref{fig:bumpontail_osde} and fig.~\ref{fig:landau_osde} are as expected such that we can proceed with the PIC simulation.

\FloatBarrier
\subsection{Density Estimation by Bilinear Interpolation}
\begin{figure}
\centering
\begin{subfigure}{0.45\textwidth}
 \includegraphics[width=\textwidth]{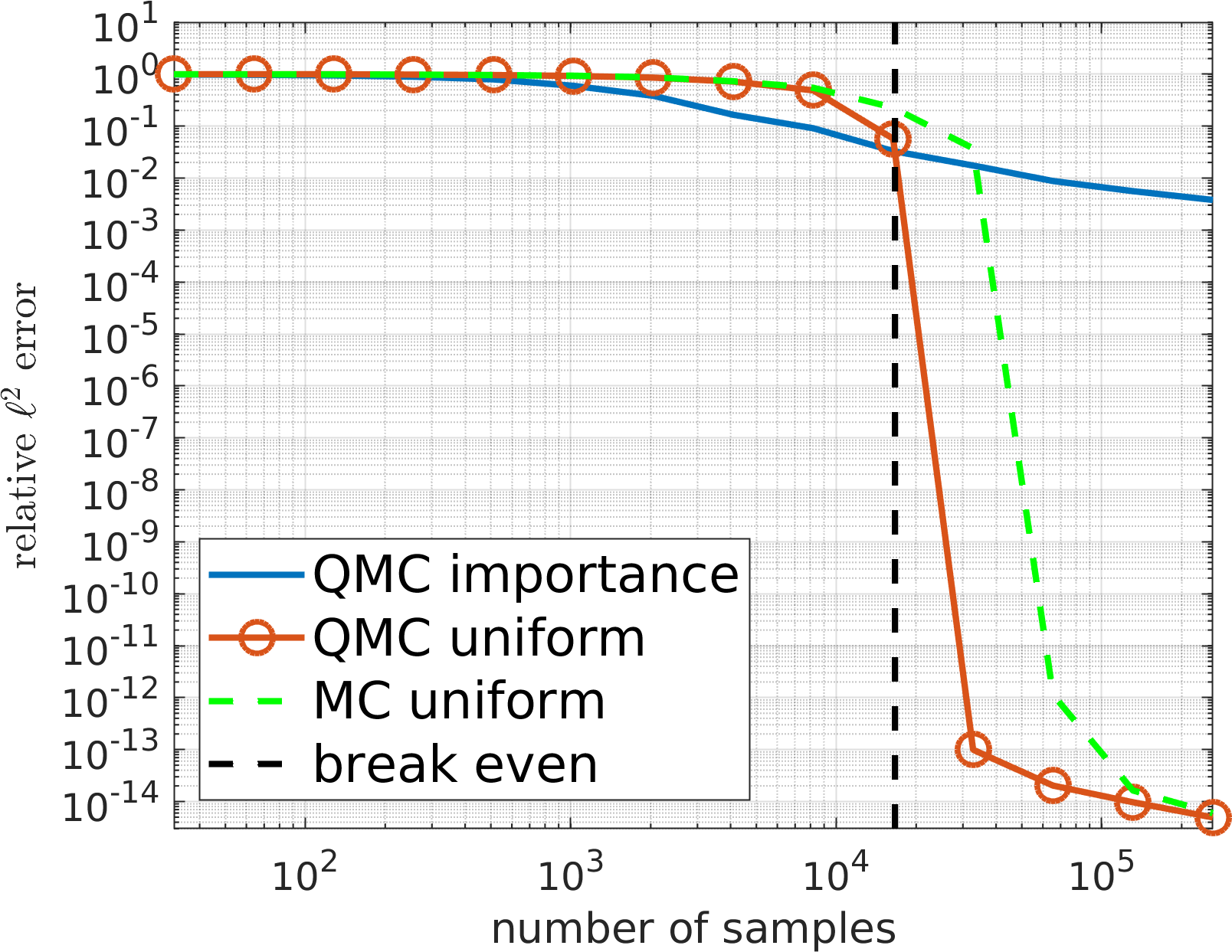} 
 \caption{Convergence of reconstruction}
\end{subfigure}
\begin{subfigure}{0.45\textwidth}
\includegraphics[width=\textwidth]{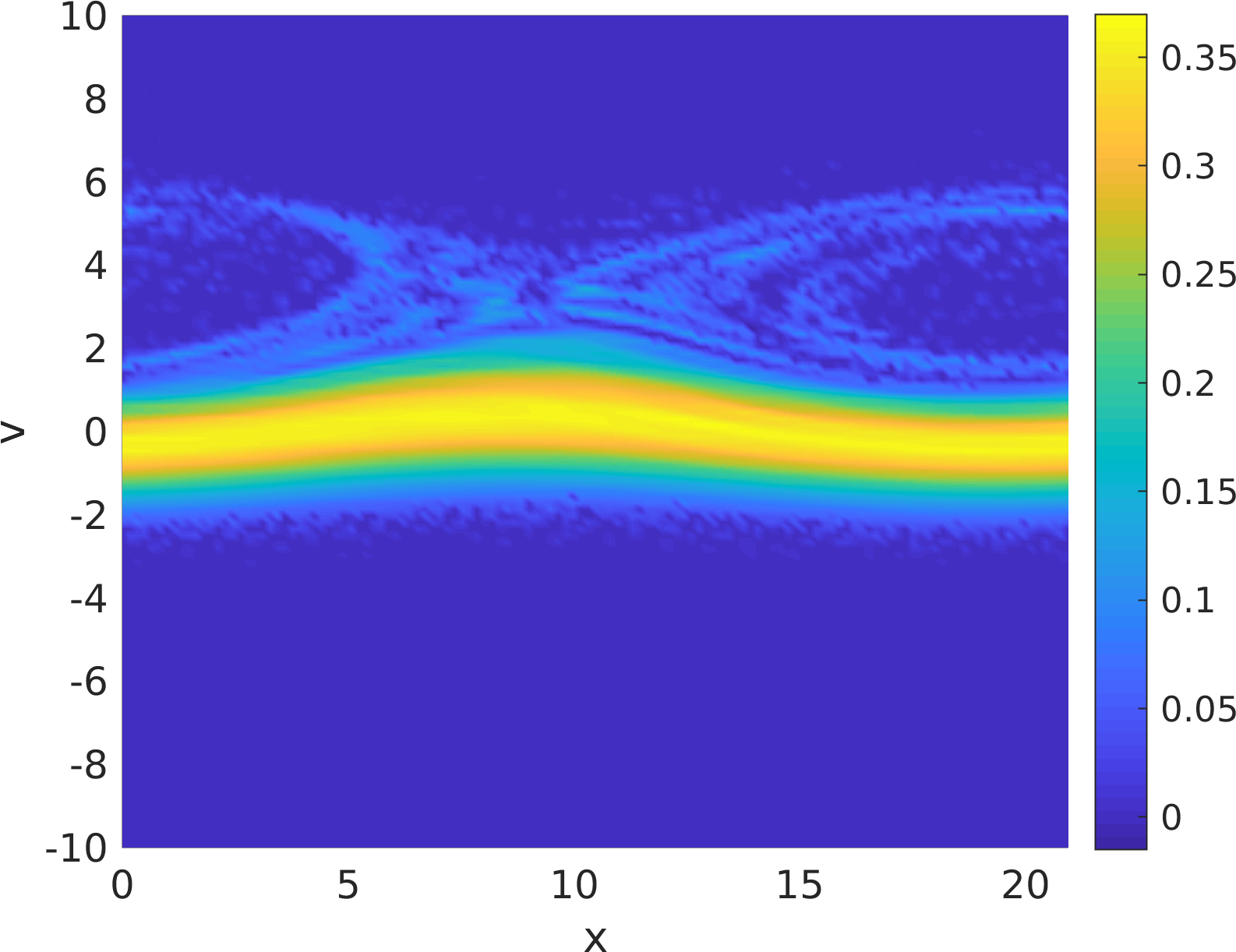}
 \caption{Bilinear reconstruction of $f$}
\end{subfigure}\\
\begin{subfigure}{0.45\textwidth}
\includegraphics[width=\textwidth]{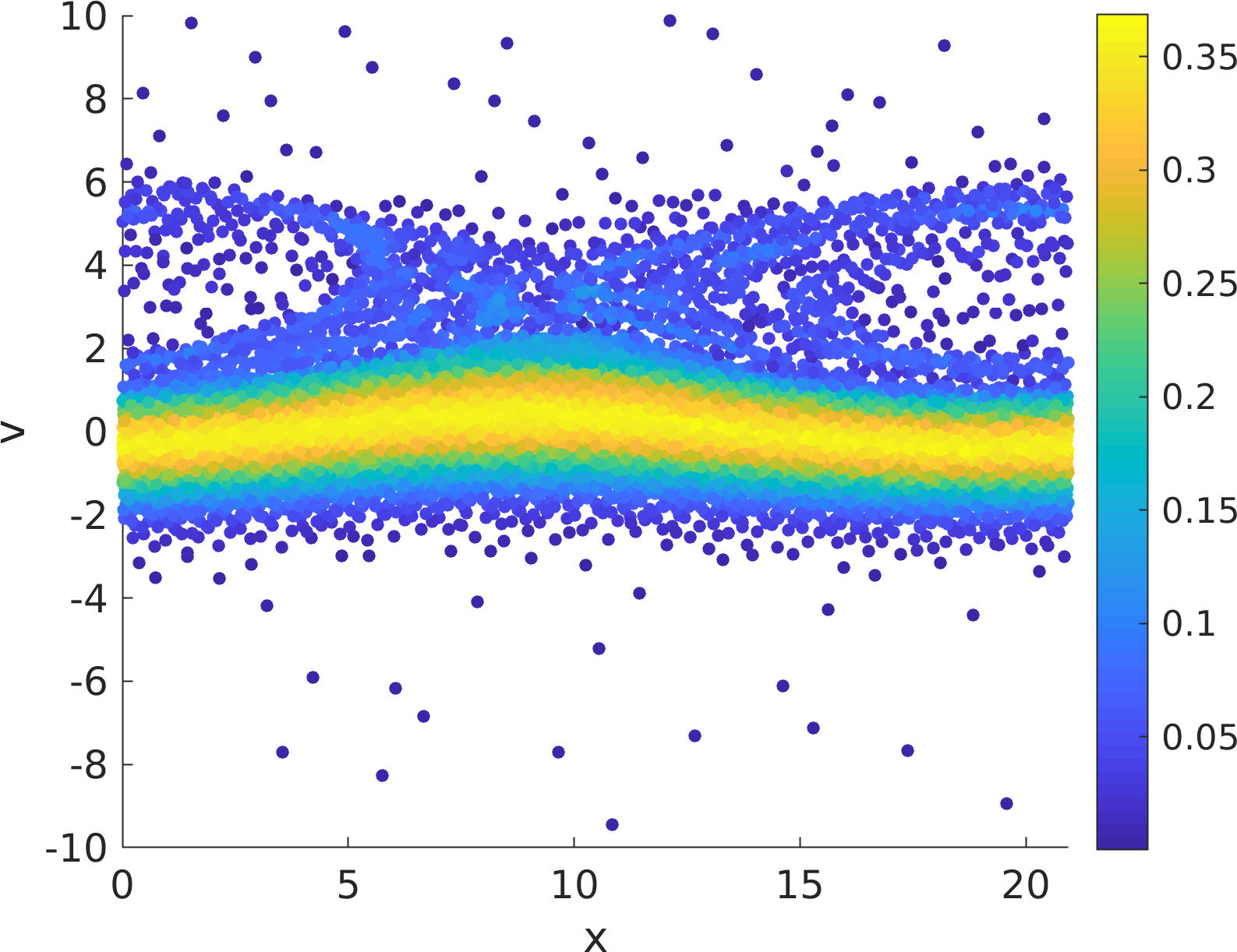}
\caption{Quasi-Random samples}
\end{subfigure}
\begin{subfigure}{0.45\textwidth}
\includegraphics[width=\textwidth]{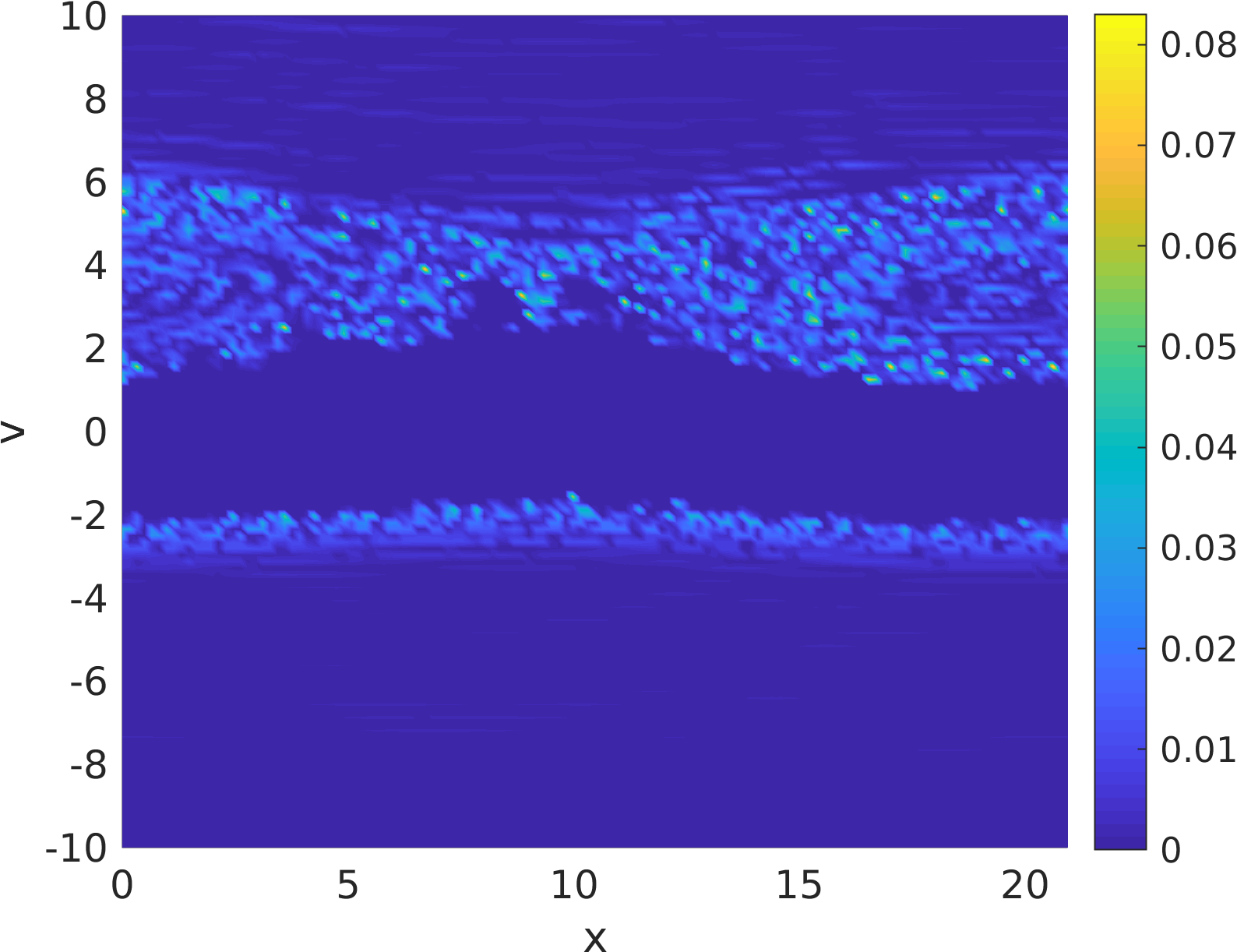}
\caption{absolute error on  of $f$}
\end{subfigure}

\caption{The density of the Bump-on-tail instability calculated by a spectral solver $N_x=N_v=128,~ \Delta t=0.01$ at $t=100$ is sampled with $N_p=10^4$ quasi-random particles in (c). Figure (b) shows the reconstruction
from these markers by bilinear interpolation with ridge regression. Due to the importance sampling
the error (d) is dominant in the regions of low but nonzero density. In (a) the error of the reconstruction by interpolation is shown for a different number of samples, where the break-even point lies at the number of degrees of freedom $N_x N_v$. When there are fewer markers than grid points and the problem is ill-posed importance sampling performs best but it does not converge for the well-posed problem. Uniform sampling is only better in the well-posed regime, where it converges to machine precision. It becomes also clear that the Quasi-Random points are more uniform than the random ones.}
\label{fig:bumpontail_interp}
\end{figure}

Using a Monte Carlo based density estimator for recovering the density $f$ from a distribution of markers is only using the information $\frac{f}{g}$ and neglecting the fact that each marker
already transports the value of the density $f$ or $g$ respectively. Using this additional information leads
to an interpolation problem. We sampled from a bilinear interpolant such that it is reasonable to test the least square fit of the bilinear interpolation coefficients under different marker distributions. 
If there are less more grid points than markers the problem is well-posed, but otherwise one has to add a 
regularization to the least square problem. The easiest choice was $L^2$ regularization known as ridge regression. Figure~\ref{fig:bumpontail_interp} shows as expected that interpolation delivers better results than Monte Carlo OSDE. For a high precision reconstruction interpolation relies on the uniformity of the interpolation points, such that uniform sampling outperforms the importance sampling by magnitudes.
{h}

\FloatBarrier
\subsection{Eulerian and Lagrangian Coupling}
\begin{figure}
\centering
\begin{subfigure}{0.7\textwidth}
\includegraphics[width=\textwidth]{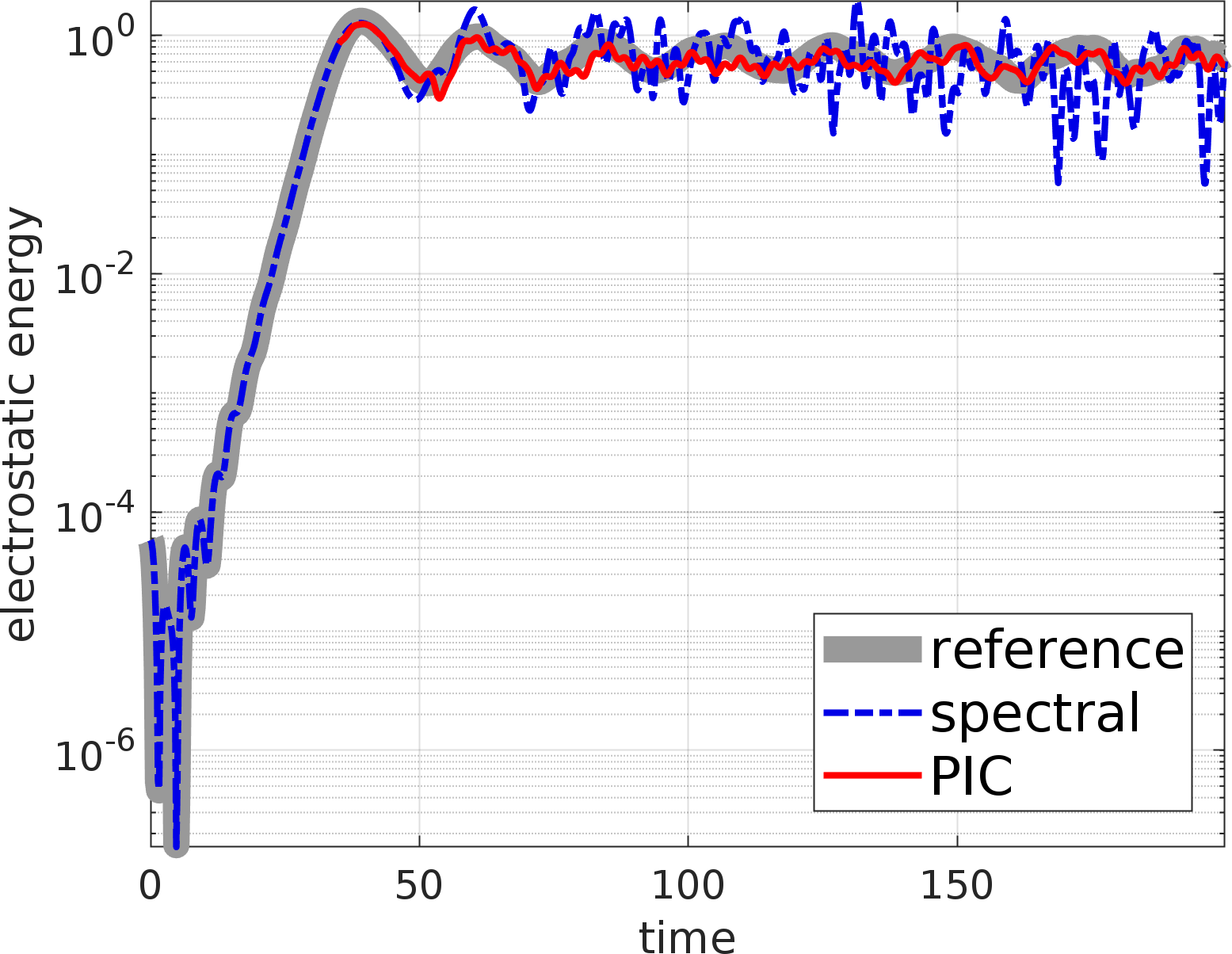}
\caption{electrostatic energy}
\label{fig:coupling:bumpontail:fieldenergy}
\end{subfigure}\\
\begin{subfigure}{0.45\textwidth}
 \includegraphics[width=\textwidth]{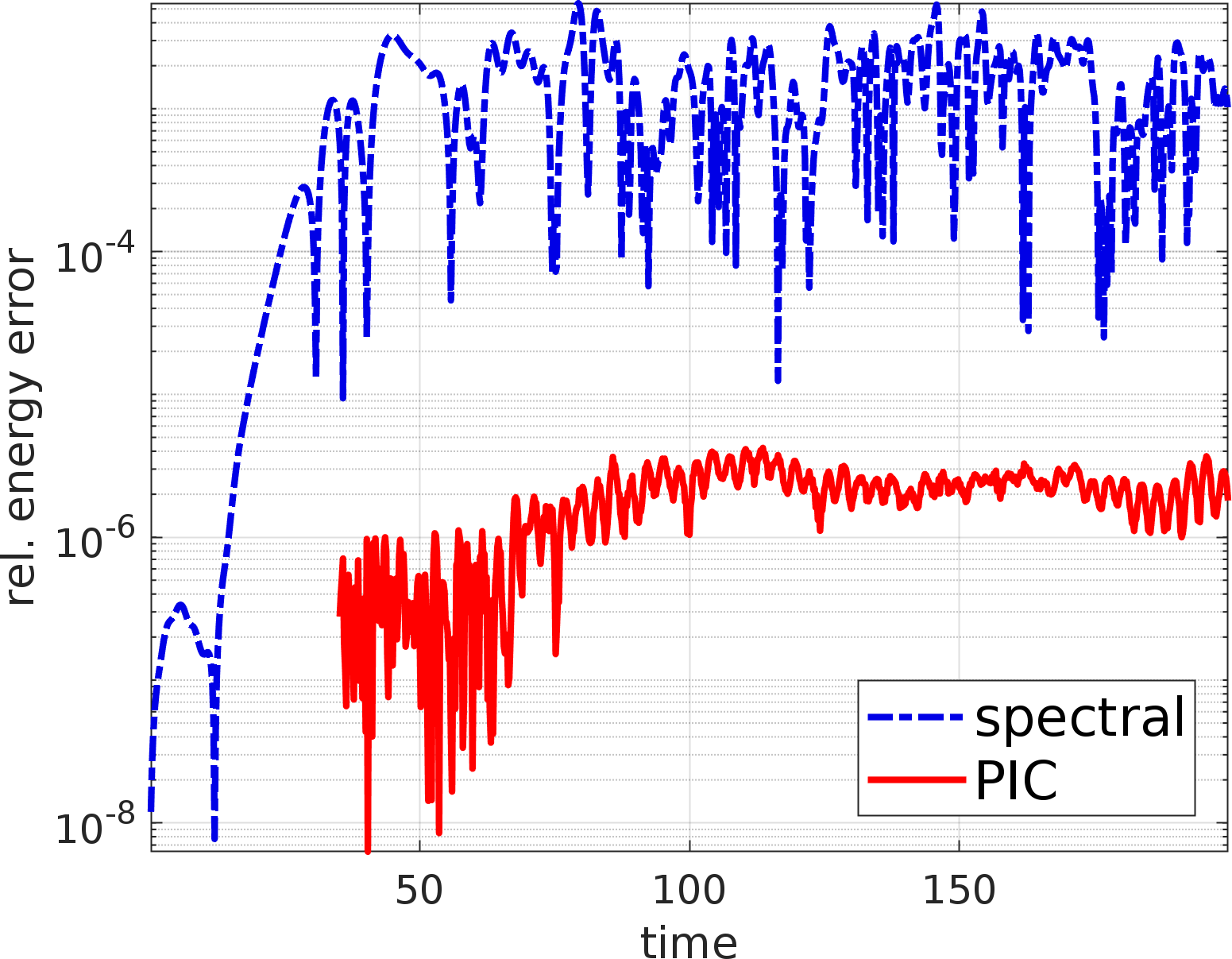} 
 \caption{relative energy error}
 \label{fig:coupling:bumpontail:energy_error}
\end{subfigure}
\begin{subfigure}{0.45\textwidth}
\includegraphics[width=\textwidth]{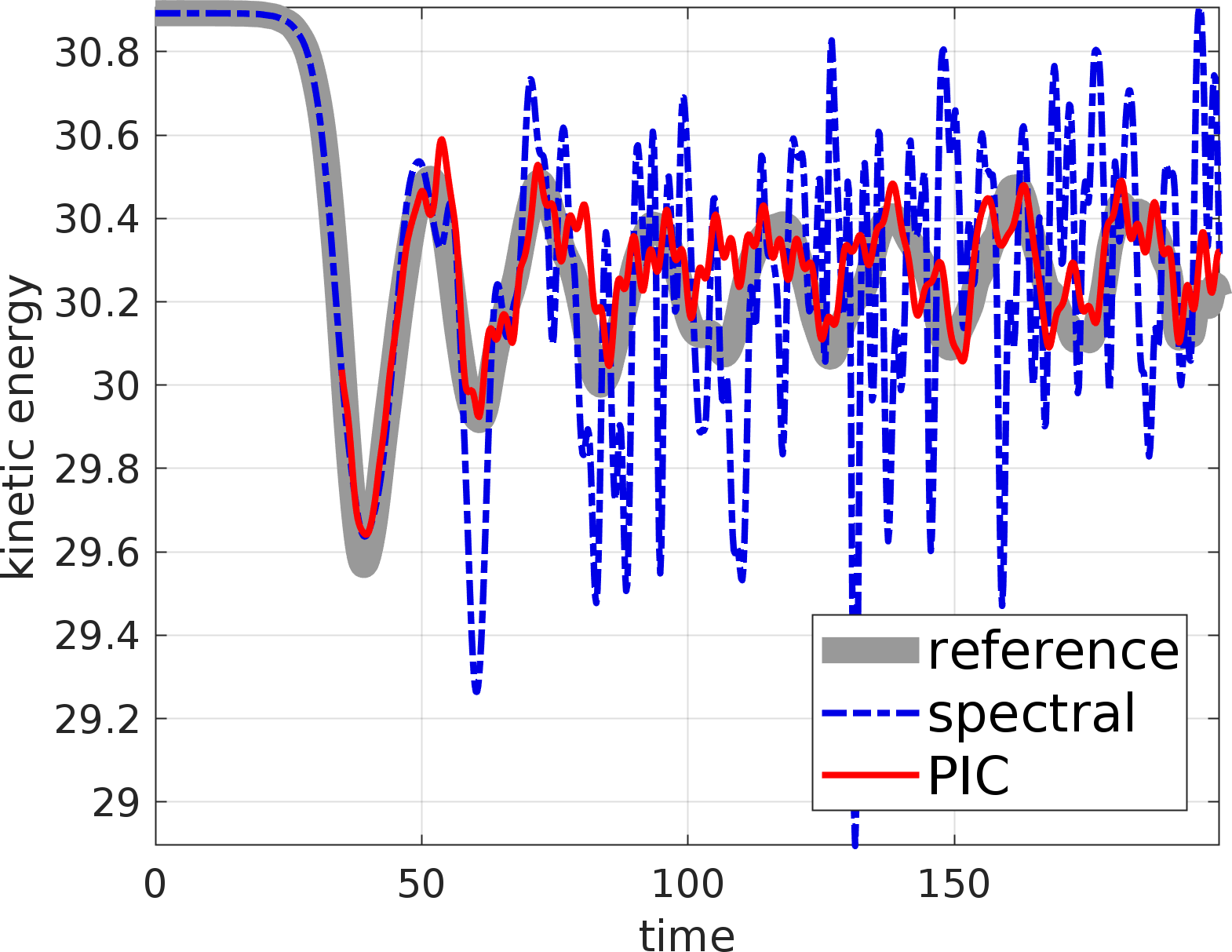}
\caption{kinetic energy}
\label{fig:coupling:bumpontail:kineticenergy}
\end{subfigure}
\caption{Transition from a spectral solver to PIC for the Bump-on-tail instability. The PIC code is initialized at $t_0=35$ by QMC inverse transform sampling.
         Although it appears from (a) and (b) as if the spectral solver is noisier than the PIC code this effect caused by filamentation is suppressed by better resolution in the reference.
         Depending on the implementation, the spectral solver, in general, outperforms the PIC code in two dimensions $(32^2<512^2<10^6)$.}
\label{fig:coupling:bumpontail}
\end{figure}

\begin{figure}
\centering
 \includegraphics[width=0.7\textwidth]{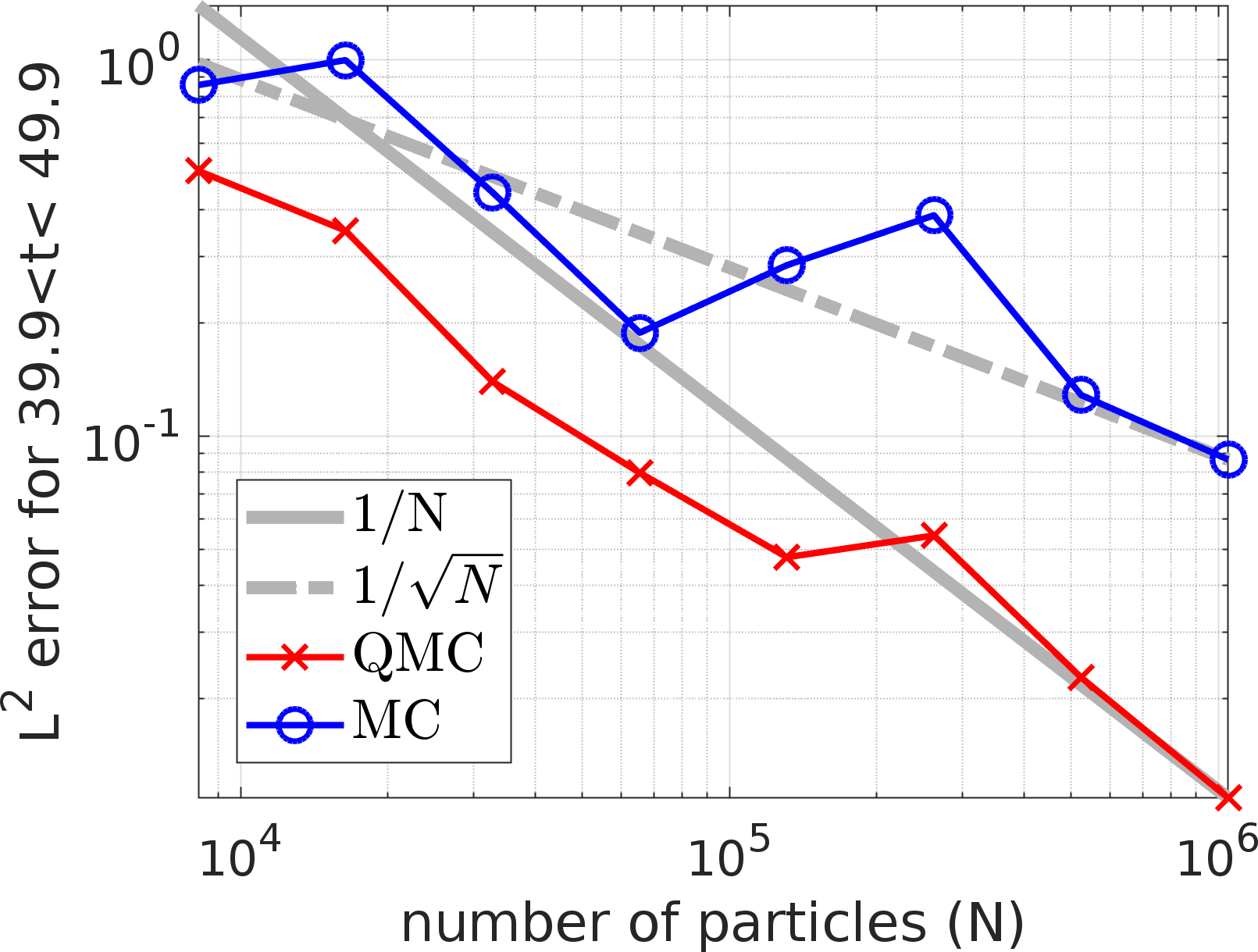} -
\caption{$L^2$ error of the electrostatic energy for quasi random (QMC) and random (MC) numbers in the spectral-PIC coupling for bump-on-tail instability from fig.~\ref{fig:coupling:bumpontail}.
         In order to dominantly include nonlinear effects the error is only taken for $t \in (39.9, 49.9)$. Here QMC performs better even with the included nonlinear effects and the excluded initial sampling.}
 \label{fig:coupling:convgc:bumpontail}
\end{figure}

\begin{figure}
\centering
\begin{subfigure}{0.45\textwidth}
 \includegraphics[width=\textwidth]{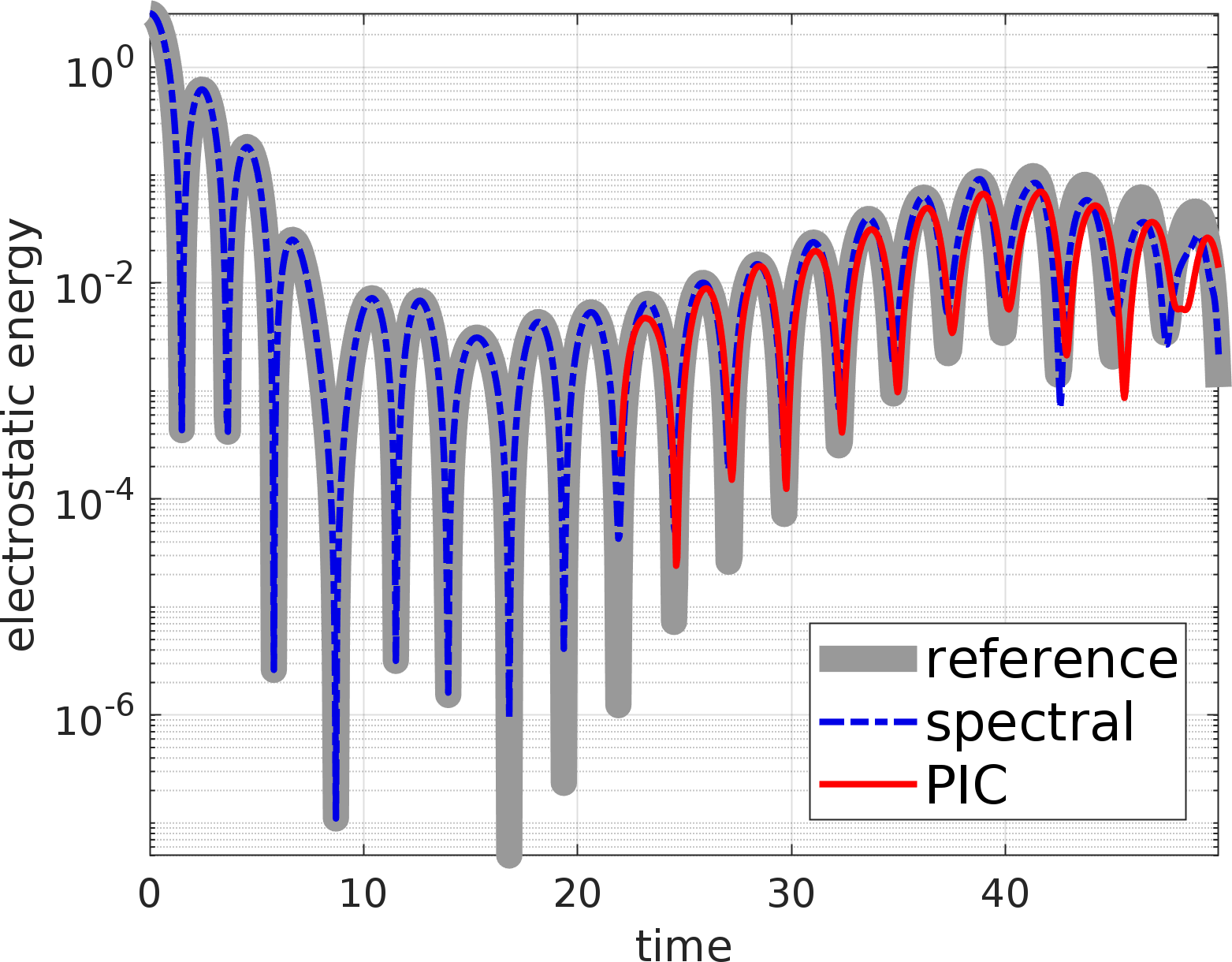} 
\caption{electrostatic energy}
\end{subfigure}
\begin{subfigure}{0.45\textwidth}
\includegraphics[width=\textwidth]{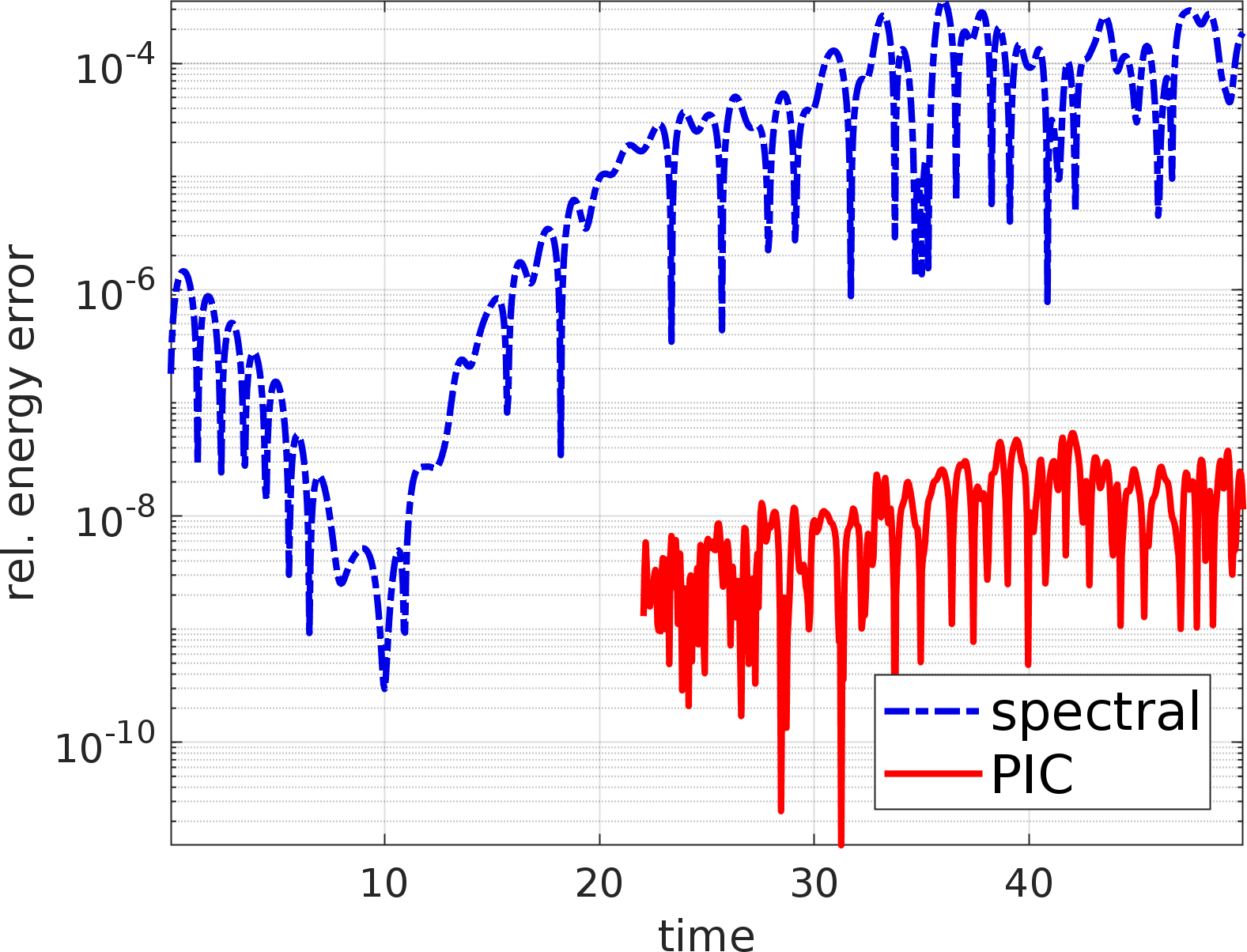}
 \caption{relative energy error}
 \end{subfigure}
\caption{Transition from a spectral solver to PIC for nonlinear Landau damping. The transition at $t_0=30$ is precisely chosen in moment of small amplitudes in order to raise the difficulty for the PIC.}
\label{fig:coupling:landau}
\end{figure}

Since we are confident that particles can be sampled correctly, we proceed with the coupling between the Eulerian and Lagrangian solver. Coupling Eulerian and PIC codes is a rapidly developing field~\cite{dominski2018tight,choi2018coupling}.
First the Bump-on-tail instability is considered, which starts with a very small amplitude posing no problem for the spectral
solver $(N_x=N_v=32,~v_{\max}=-v_{\min}=10,~\Delta t=0.1)$, which follows the reference solution $(N_x=N_v=8 \times 32= 512)$ properly in the linear phase~\ref{fig:coupling:bumpontail:fieldenergy}
but later suffers from oscillations due to filamentations and aliasing. Therefore, we switch to PIC ($N_p=10^6~,N_f=16$) at $t_0=35$.
Since the resolution is quite low and we want to suppress aliasing due to the low order interpolation the spectral density is zero padded with a factor of $N_{\mathrm{pad}}=32$ yielding
a $1024 \times 1024$ phase space grid for the sampler. Given the curse of dimensionality high order splines would be a better solution.
Nevertheless the PIC code follows the reference visibly better in fig.~\ref{fig:coupling:bumpontail:kineticenergy} and fig.~\ref{fig:coupling:bumpontail:kineticenergy}.
The third order symplectic Runge Kutta time discretization is the same for PIC and the spectral solver, but
PIC is derived from a Lagrangian formulation such that the better energy conservation in fig.~\ref{fig:coupling:bumpontail:energy_error} is no surprise.\\
Although fig.~\ref{fig:bumpontail_osde} and fig.~\ref{fig:landau_osde} confirm the superior convergence rate of the QMC numbers for the initial sampling, the question remains whether this stays true
over the nonlinear phase. This question was already addressed before~\cite{ameres2016adv,ameres2018} and can also be answered positively here by fig.~\ref{fig:coupling:convgc:bumpontail}.
Note that the symplectic Runge Kutta scheme preserves phase space volume, hence the Jacobian of the discrete flux is exactly one which preserves also the likelihood of each marker.
This means that the discrete flux induces a measure-preserving map such that the Hausdorff measure is preserved. Therefore, the uniformity of the low discrepancy sequence is preserved such that
the higher order convergence rates for QMC keep their validity~\cite{basu2016transformations,aistleitner2013low,de2018quasi}.
The only confusing issue that can emerge is, that the total variation (QMC) as well as the variance (MC) 
of the \textit{entire} map from the initial condition to a certain time $t$ increases with the nonlinearities of the transport. This, however, does not change the convergence rates at a certain time $t$.\\
The same procedure with $(N_x=N_v=64,~v_{\max}=-v_{\min}=6.5,~\Delta t=0.05,~t_0=30,~N_p=10^6,~N_f=16, N_{\mathrm{pad}}=32$ is repeated for nonlinear Landau damping.
Here, because of the many perturbations at small amplitudes the difference between the spectral solver and PIC remains small, see fig.~\ref{fig:coupling:landau}.

\FloatBarrier
\section{Conclusion and Outlook}
\begin{figure}[h]
\centering
\includegraphics[width=0.75\textwidth]{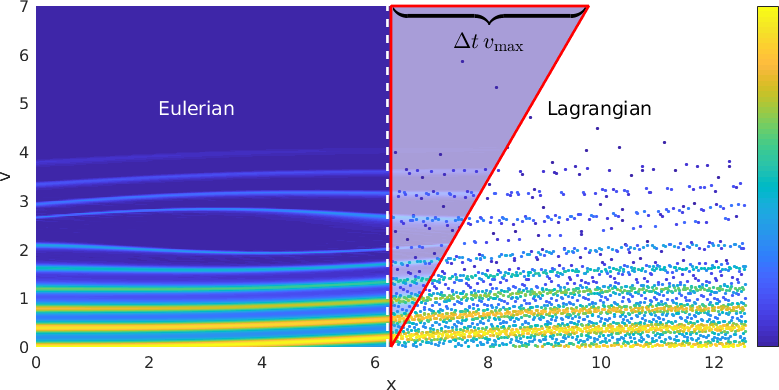}
 \caption{The image shows the upper half of the phase space of a Vlasov--Poisson simulation with a transition from a Eulerian to a Lagrangian region. This means the entire phase-space is flowing to the right. The white dashed line marks the boundary/interface position. During the split-step $\dot{X}=V$ the red triangle encloses the area of phase space entering the Lagrangian domain at each time step of size $\Delta t$. In this case, the Eulerian solver can also precalculate the volume in the red triangle, such that it can be sampled directly via inverse transform sampling at each time step, which results in a smooth density. If the boundary condition is only given at a fixed point $\bar{x}$ with density
 $f_{\mathbf{X}=\bar{x}}(v)$ it is not sufficient to draw the particles velocity as $\mathbf{V} \sim f_{\mathbf{X}=\bar{x}} $ and set $\mathbf{X}=\bar{x}$. To account for a volume one has to draw 
 $\mathbf{X} \sim \mathcal{U}(\bar{x}, \bar{x} + \Delta t \mathbf{X})  $ uniformly. For a strong electric field, the splitting might not be appropriate such that the phase space integral becomes more complicated but the principle stays the same.}
\label{fig:eulerianlagrangianboundary}
\end{figure}
It was shown how to use Quasi Monte Carlo numbers to sample from arbitrary phase-space densities. Not only does this improve convergence but opens up the possibility to couple established Eulerian codes to the big PIC codes for kinetic and gyrokinetic systems.
It should also be noted that 
Considering that Vlasov problems are high dimensional one has to realize that computing high dimensional marginals can be expensive. 
Here~\cite{olver2013fast} proposed to use a low-rank spectral representation, which is a major underlying idea of Approxfun.jl~\cite{olver2014practical}.
For the Vlasov--Poisson system there exist solvers in tensor train format~\cite{kormann2015solving} that already provide a low-rank approximation of $f$ 
which possibly can be exploited for efficient sampling in the future.
Another aspect not followed here is the transition from a Lagrangian to a Eulerian solver. OSDE is for sure a possibility but this misses the point of Lagrangian solvers completely such that an interpolation approach as shown in fig.~\ref{fig:bumpontail_interp} is recommended.
Each particle transports the value of the density along the characteristics such that the key for such a transition lies in a suitable interpolation, which has already been extensively discussed for Semi-Lagrangian solvers~\cite{sonnendrucker1999semi,bermejo1991analysis}.
In two dimensions phase-space conservation is the same as symplecticity, but in higher dimensions, symplecticity is something mildly stronger. Hence, it remains the question, what impact symplecticity has on the discrepancy.\\
Another important application of the sampling techniques presented here is the implementation of boundary conditions for the Vlasov equation in PIC. The standard approach is to draw particles at an interface position according to a velocity distribution at each time step.  PIC relies on phase space conservation such that any insertion of markers should actually be an insertion of a volume of phase space. This means that the time step cannot be assumed to be infinitesimally small, but also the boundary condition has to be integrated exactly over time. In most cases, this can be done analytically and otherwise, a numerical pre-calculation is sufficient. As explained in fig.~\ref{fig:eulerianlagrangianboundary} this requires sampling nontrivial phase space volumes, where this article opens new possibilities.

\section{Acknowledgement}
This work has been carried out within the framework of the EUROfusion Consortium and has received funding from the Euratom research and training programme 2014-2018 and 2019-2020 under grant agreement No 633053. The views and opinions expressed herein do not necessarily reflect those of the European Commission.

\FloatBarrier

\printbibliography

\end{document}